\documentclass[10pt,twocolumn,twoside]{IEEEtran}
\usepackage{times} 
\usepackage{pdfsync} 
\usepackage{epsfig}
\usepackage{epstopdf}
\usepackage{enumerate}
\usepackage[cmex10]{amsmath} \interdisplaylinepenalty=2500
\usepackage{amssymb}
\usepackage{color} 
\usepackage{cite} 
\usepackage{array}
\usepackage{mathrsfs}
\usepackage{fix2col} 
\usepackage{hyperref}
\usepackage{subfig}

\graphicspath{{Figures/}}

\newcommand{\Prob}{\mathbb{P}}
\newcommand{\E}{\mathbb{E}}

\DeclareMathOperator{\tr}{tr}

\def\qed
{\hfill\vbox{\hrule width 0.5em\nointerlineskip\hbox to
0.5em{\vrule height 0.5em \hfill\vrule height
0.5em}\nointerlineskip\hrule width 0.5em}}

\newtheorem{thm}{Theorem}[section]

\newtheorem{lem}{Lemma}[section]

\newtheorem{remark}{Remark}[section] 

\title{Optimal Energy Allocation for Kalman Filtering over Packet Dropping Links with Imperfect Acknowledgments and Energy Harvesting Constraints\thanks{A preliminary version of this paper was presented at the 4th IFAC NecSys workshop, Koblenz, Germany, Sep. 2013.}}

\author{Mojtaba Nourian, Alex S. Leong and Subhrakanti Dey 
\thanks{M. Nourian and A.S. Leong are with the Department of Electrical and Electronic Engineering, The University of Melbourne, VIC 3010, Australia (email: {\tt\small \{mojtaba.nourian,asleong\}@unimelb.edu.au}). 

S. Dey is with the Dept of Eng. Sciences in Uppsala University, Sweden 
(email: {\tt\small Subhra.Dey@signal.uu.se}).}
}

\begin{document}

\maketitle

\begin{abstract} 
This paper presents a design methodology for optimal transmission energy allocation at a sensor equipped with energy harvesting technology for remote state estimation of linear stochastic dynamical systems.
In this framework, the sensor measurements as noisy versions of the system states are sent to the receiver over a packet dropping communication channel. The packet dropout probabilities of the channel depend on both the sensor's transmission energies and time varying wireless fading channel gains. The sensor has access to an energy harvesting source which is an everlasting but unreliable energy source compared to conventional batteries with fixed energy storages. The receiver performs optimal state estimation with random packet dropouts to minimize the estimation error covariances based on received measurements. The receiver also sends packet receipt acknowledgments to the sensor via an erroneous feedback communication channel which is itself packet dropping. 

The objective is to design optimal transmission energy allocation at the energy harvesting sensor to minimize either a finite-time horizon sum or a long term average (infinite-time horizon) of the trace of the expected estimation error covariance of the receiver's Kalman filter. These problems are formulated as Markov decision processes with imperfect state information. 
The optimal transmission energy allocation policies are obtained by the use of dynamic programming techniques. Using the concept of submodularity, the structure of the optimal transmission energy policies are studied. Suboptimal solutions are also discussed which are far less computationally intensive than optimal solutions.
Numerical simulation results are presented illustrating the performance of the energy allocation algorithms.
\end{abstract}


\section{Introduction}
\IEEEPARstart{W}{ireless} sensor network (WSN) technologies arise in a wide range of applications such as environmental data gathering \cite{akyildiz2002survey,duarte2003data}, mobile robots and autonomous vehicles \cite{lamarca2002making,Chong2003}, and monitoring of smart electricity grids \cite{gungor2010opportunities,erol2011wireless}, among many others. In these applications one of the important challenges is to improve system performance and reliability under resource (e.g., energy/power, computation and communication) constraints.

A considerable amount of research has recently been devoted to the concept of energy harvesting \cite{RoundyWrightRabaey} (see also \cite{niyato2007wireless,sharma2010optimal,Ozel11,ho2012optimal,kashef2012optimal,NayyarBasar} among other papers). This is motivated by energy limited WSN applications where sensors may need to operate continuously for years on a single battery. In the energy harvesting paradigm the sensors can recharge their batteries by collecting energy from the environment, e.g. solar, wind, water, thermal or mechanical vibrations. However, the amount of energy harvested is random as most renewable energy sources are unreliable. In this work we will consider the remote Kalman filtering problem with random packet dropouts and imperfect receipt acknowledgments when the sensors are equipped with energy harvesting technology, and as a result, are subject to energy harvesting constraints.

Since the seminal work of \cite{Sinopoli}, the problem of state estimation or Kalman filtering over packet dropping communication channels has been studied extensively (see for example\cite{LiuGoldsmith,XuHespanha,HuangDey,Epstein_Automatica,Schenato,MoSinopoli,Quevedo_Automatica} among others). The reader is also referred to the comprehensive survey \cite{Schenato_proceedings} for some of the research on the area of control and estimation over lossy networks up to 2007. In these problems sensor measurements (or state estimates in the case of \cite{XuHespanha}) are grouped into packets which are transmitted over a packet dropping link such that either the entire packet is received or lost in a random manner. The focus in these works is on deriving conditions on the packet arrival rate in order to guarantee the stability of the Kalman filter. 

There are other works which are concerned with estimation performance (e.g. minimizing the expected estimation error covariance) rather than just stability. For instance, power allocation techniques\footnote{We measure energy on a per channel use basis and we will refer to energy and power interchangeably.} (without energy harvesting constraints) have been applied to the Kalman filtering problem in \cite{QuevedoAhlenOstergaard,Alex_CDC12,ShiChengChen} in order to improve the estimation performance. In these works energy allocation can be used to improve system performance and reliability. 

In conventional wireless communication systems, the sensors have access either to a fixed energy supply or have batteries that may be easily rechargeable/replaceable. Therefore, the sum of energy/power constraint is used to model the energy limitations of the battery-powered devices (see \cite{Alex_CDC12}). However, in the context of WSNs the use of energy harvesting is more practical, e.g., in remote locations with restricted access to an energy supply, and even essential where it is dangerous or impossible to change the batteries \cite{ho2010markovian,ho2012optimal}. In these situations it is possible to have communication devices with on-board energy harvesting capability which may recharge their batteries by collecting energy from the environment including solar, thermal or mechanical vibrations.

Typically, the harvested energy is stored in an energy storage such as a rechargeable battery which then is used for communications or other processing. Even though the energy harvesters provide an everlasting energy source for the communication devices, the amount of energy expenditure at every time slot is constrained by the amount of stored energy currently available. This is unlike the conventional communication devices that are subject only to a sum energy constraint. Therefore, a causality constraint is imposed on the use of the harvested energy \cite{ho2012optimal}. 
Communication schemes for optimizing throughput for transmitters with energy harvesting capability have been studied in \cite{Ozel11,ho2012optimal}, while a remote estimation problem with an energy harvesting sensor was considered in \cite{NayyarBasar} which minimized a cost consisting of both the distortion and number of sensor transmissions.

In this paper we study the problem of optimal transmission energy allocation at an energy harvesting sensor for remote state estimation of linear stochastic dynamical systems. In this model, sensor's measurements as noisy versions of the system's states are sent to the receiver over a packet dropping communication channel. Similar to the channel models in \cite{Quevedo_Automatica}, the packet dropout probabilities depend on both the sensor's transmission energies and time varying wireless fading channel gains. The sensor has access to an energy harvesting source which is an everlasting but unreliable energy source compared to conventional batteries with fixed energy storages. The receiver performs a Kalman filtering optimal state estimation with random packet dropouts to minimize the estimation error covariances based on received measurements. In general, knowledge at the sensor of whether its transmissions have been received at the receiver is usually achieved via some feedback mechanism. Here, in contrast to the models in \cite{Alex_CDC12,Nourian_Necsys13} the feedback channel from receiver to sensor is also a packet dropping erroneous channel leading to a more realistic formulation. The energy consumed in transmission of a packet 
is assumed to be much larger than that for sensing or processing at the sensor and thus energy consumed in sensing and processing is not taken into account 
in our formulation.

The objective of this work is to design optimal transmission energy allocation (per packet) at the energy harvesting sensor to minimize either a finite-time horizon sum or a long term average (infinite-time horizon) of the trace of the expected estimation error covariance of the receiver's Kalman filter. The important issue in this problem formulation is to address the trade-off between the use of available stored energy to improve the current transmission reliability and thus state estimation accuracy, or storing of energy for future transmissions which may be affected by higher packet loss probabilities due to severe fading. 

These optimization problems are formulated as Markov decision processes with imperfect state information. The optimal transmission energy allocation policies are obtained by the use of dynamic programming techniques. Using the concept of submodularity \cite{topkis2001supermodularity}, the structure of the optimal transmission energy policies are studied. Suboptimal solutions which are far less computationally intensive than optimal solutions are also discussed. Numerical simulation results are presented illustrating the performance of the energy allocation algorithms.  

Previous presentation of the model considered in this paper includes \cite{Nourian_Necsys13} which investigates the case with perfect acknowledgments at the sensor. Here, we address the more difficult problem where the feedback channel from receiver to sensor is an imperfect erroneous channel modelled as an erasure channel with errors.

In summary, the main contributions of this paper are as follows:
\begin{enumerate}[i)]
\item Unlike a large number of papers focusing on the stability for Kalman filtering with packet loss, e.g. \cite{LiuGoldsmith,XuHespanha,HuangDey,Epstein_Automatica,Schenato,MoSinopoli,Quevedo_Automatica},
we focus on the somewhat neglected issue of estimation error performance (noting that stability only guarantees bounded estimation error) in the presence of packet loss and how to optimize it via power/energy allocation at the sensor transmitter. Note that it is quite common to study optimal power allocation in the context of a random stationary source estimation in {\em fading} wireless sensor networks \cite{Cui_TSP}, but this issue has received much less attention in the context of Kalman filtering over packet dropping links which are randomly time-varying. In particular, we consider minimization of a long-term average of error covariance minimization for the Kalman filter by optimally allocating energy for individual packet transmissions over packet dropping links with randomly varying packet loss probability due to fading. While a version of this problem was considered in our earlier conference paper \cite{Alex_CDC12}, we extend the problem setting and the analysis along multiple directions as described below. 

\item Unlike \cite{Alex_CDC12}, we consider an energy harvesting sensor that is not constrained by a fixed initial battery energy, but rather the randomness of the harvested energy pattern. Energy harvesting is a promising solution to the important problem of energy management in wireless sensor networks. Furthermore, recent advances in hardware have made energy harvesting technology a practical reality \cite{RoundyWrightRabaey}. 

\item We provide a new sufficient stability condition for bounded long term average estimation error, which depends on the packet loss probability (which is a function of the channel gain, harvested energy and the maximum battery storage capacity) and the statistics of the channel gain and harvested energy process. Although difficult to verify in general, we provide simpler forms of this condition in  when the channel gains and harvested energy processes follow familiar statistical models such as independent and identically distributed processes or finite state Markov chains.

\item We consider the case of imperfect feedback acknowledgements, which is more realistic but more difficult to study than the case of perfect feedback acknowledgements. We model the feedback channel by a general erasure channel with errors.

\item It is well known that the optimal solution obtained by a stationary control policy minimizing the infinite horizon control cost is computationally prohibitive. Thus motivated, we provide structural results on the optimal energy allocation policy which lead to threshold policies which are optimal and yet very simple to implement in some practical cases, e.g. when the sensor is equipped with binary transmission energy levels. Note that most sensors usually have a finite number of transmission energy/power levels and for simplicity, sensors can be programmed to only have two levels. 
 
\item Finally, also motivated by the computational burden for the optimal control solution in the general case of imperfect acknowledgments, we provide a sub-optimal solution based on an estimate of the error covariance at the receiver. Numerical results are presented to illustrate the performance gaps between the optimal and sub-optimal solutions.

\end{enumerate}

The organization of the paper is as follows. The system model is given in Section \ref{model_sec}. The optimal energy allocation problems subject to energy harvesting constraints are formulated in Section \ref{optim_sec}. In Section \ref{optimSol_sec} the optimal transmission energy allocation policies are derived by the use of dynamic programming techniques. Section \ref{Sub_sec} presents suboptimal policies which are less computationally demanding. The structure of the optimal transmission energy allocation policies are studied in Section \ref{Str_sec}. Section \ref{numerical_sec} presents the numerical simulation results. Finally, concluding remarks are stated in Section \ref{conc_sec}.

\begin{figure}[!t]
\centering 
\includegraphics[width=9.2cm,height=4cm]{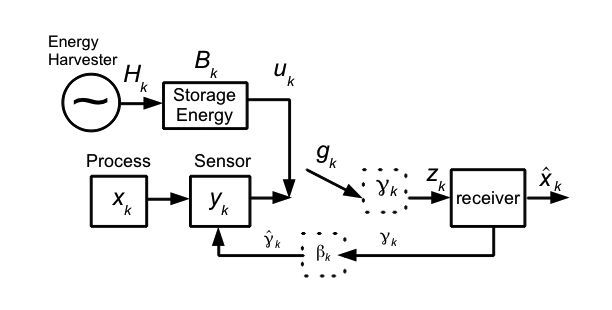} 
\caption{System model}
\label{system_model}
\end{figure}

\section{System Model} \label{model_sec}

A diagram of the system architecture is shown in Fig. \ref{system_model}. The description of each part of the system is given in detail below. 

\subsection{Process Dynamics and Sensor Measurements} \label{ProcessDyn_sec}
We consider a linear time-invariant stochastic dynamical process
\begin{align*}
& x_{k+1} = A x_k+w_k, \qquad k \geq 0
\end{align*}
where $x_k \in \mathbb{R}^n$ is the process state at time $k \geq 0$, $A \in \mathbb{R}^{n \times n}$, and $\{w_k, k \geq 0\}$ is a sequence of independent and identically distributed (i.i.d.) Gaussian noises with zero mean and positive definite covariance matrix $Q >0$. The initial state of the process $x_0$ is a Gaussian random vector, independent of the process noise sequence $\{w_k, k \geq 0\}$, with mean $\bar x_0 := \E x_0$ and covariance matrix $P_{x_0}$.

The sensor measurements are obtained in the form
\begin{align*}
& y_k=C x_k+v_k, \qquad k \geq 0
\end{align*}
where $y_k \in \mathbb{R}^{m}$ is the observation at time $k \geq 0$, $C \in \mathbb{R}^{m \times n}$, and $\{v_k, k \geq 0\}$ is a sequence of i.i.d. Gaussian noises, independent of both the initial state $x_0$ and the process noise sequence $\{w_k, k \geq 0\}$, with zero mean and a positive semi-definite covariance matrix $R \geq 0$.

We enunciate the following assumption:

({\bf A1}) We assume that $(A,Q^{1/2})$ is stabilizable and $(A,C)$ is detectable.  \qed

\subsection{Forward Communication Channel} \label{FCchannel_sec}

The measurement $y_k$ is then sent to a receiver over a packet dropping communication channel such that $y_k$ (considered as a packet) is either exactly received or the packet gets lost due to corrupted data or substantial delay. The packet dropping channel is modelled by
\begin{align*}
& z_k=\gamma_k y_k, \qquad k \geq 0
\end{align*}
where $z_k$ is the observation obtained by the receiver at time $k$, and $\gamma_k=1$ denotes that the measurement packet is received, while $\gamma_k=0$ denotes that the packet containing the measurement $y_k$ is lost. 

Similar to \cite{Quevedo_Automatica}, we adopt a model for the packet loss process $\{\gamma_k\}$ that is governed by the time-varying wireless fading channel gains $\{g_k\}$ and sensor transmission energy allocation (per packet) $\{u_k\}$ over this channel. In this model, the conditional packet reception probabilities are given by
\begin{equation}
\label{h_defn}
\Prob(\gamma_k=1|g_k,u_k) \triangleq h(g_k u_k) 
\end{equation}
where $h(\cdot): [0,\infty) \rightarrow [0,1]$ is a monotonically increasing continuous function. The form of $h(\cdot)$ will depend on the particular digital modulation scheme being used \cite{Proakis}. 

We consider the case where the set of fading channel gains $\{g_k\}$ is a first-order stationary and homogeneous Markov fading process (see \cite{Quevedo_TAC13}) where the channel remains constant over a fading block (representing the coherence time of the channel \cite{Rappaport}). Note that the stationary first-order Markovian modelling includes the case of independent and identically distributed (i.i.d.) processes as a special case. 

We assume that channel state information is available at the transmitter such that it knows the values of the channel gains $g_k$ at time $k$. In practice, this can be achieved by channel reciprocity between the sensor-to-receiver and receiver-to-sensor channels (such as in typical time-division-duplex (TDD) based transmissions). In this scenario, the sensor can estimate the channel gain based on pilot signals transmitted from the remote receiver at the beginning of each fading block. Another possibility (if channel reciprocity does not hold) is to estimate the channel at the receiver based on pilot transmissions from the sensor and send it back to the sensor by channel state feedback. However, transmitting pilot signals consumes energy which should then be taken into account. To conform with our problem formulation, we therefore assume that channel reciprocity holds.

\subsection{Energy Harvester and Battery Dynamics}

Let the unpredictable energy harvesting process be denoted by $\{H_k\}$ which is also modelled as a stationary first-order homogeneous 
Markov process, and which is independent of the fading process $\{g_k\}$. This modelling for the harvested energy process is justified by empirical measurements in the case of solar energy \cite{ho2010markovian}. 

We assume that the dynamics of the stored battery energy $B_{(\cdot)}$ is given by the following first-order Markov model
\begin{align} \label{EHM}
& B_{k+1} = \min \{B_k-u_k+H_{k+1},B_{\textrm{max}}\}, \quad k \geq 0
\end{align}
with given $0 \leq B_0 \leq B_{\textrm{max}}$, where $B_{\textrm{max}}$ is the maximum stored energy in the battery.

\subsection{Kalman Filter at Receiver}

The receiver performs the optimal state estimation by the use of Kalman filtering based on the history $\mathcal F_k := \sigma(z_t,\gamma_t, 0 \leq t \leq k)$ which is the $\sigma$-field generated by the available information at the receiver up to time $k$. We use the convention $\mathcal F_0 =\{\O,\Omega\}$.

The optimal Kalman filtering and prediction estimates of the process state $x_k$ are given by $\hat{x}_{k|k} =  \mathbb{E} [x_{k} | \mathcal F_k]$ and $ \hat{x}_{k+1|k} =  \mathbb{E} [x_{k+1} | \mathcal F_k]$, respectively. The corresponding Kalman filter error covariances are defined as
\begin{align*}
& P_{k|k}  =  \mathbb{E} [(x_{k} \!-\! \hat{x}_{k|k})(x_{k} \!-\! \hat{x}_{k|k})^T |  \mathcal F_k]\\
& P_{k+1} := P_{k+1|k}  =  \mathbb{E} [(x_{k+1} - \hat{x}_{k+1|k})(x_{k+1} - \hat{x}_{k+1|k})^T | \mathcal F_k].
\end{align*}
The Kalman recursion equations for $\hat{x}_{k|k}$ and $\hat{x}_{k+1|k}$ are given in \cite{Sinopoli}. In this paper we focus on the estimation error covariance $P_k$ which satisfies the random Riccati equation
\begin{align} \label{EH:RiccatiEq}
& \! \! \!  P_{k+1}\! = \! A P_k A^T \!+\! Q \!-\! \gamma_k A P_k C^T (C P_k C^T \!+\! R)^{-1} C P_k A^T
\end{align}
for $k \geq 0$ where $P_0 = \E[(x_0-\bar x_0)(x_0- \bar x_0)^T]=P_{x_0}$ (see \cite{Sinopoli}). Note that $\gamma_k$ appears as a random coefficient in the Riccati equation (\ref{EH:RiccatiEq}). Since (i) the derivation in \cite{Sinopoli} allows for time-varying packet reception probabilities, and (ii) in the model of this paper the energy allocation $u_k$ only affects the probability of packet reception via (\ref{h_defn}) and not the system state that is being estimated, the estimation error covariance recursion is of the form (\ref{EH:RiccatiEq}) as given in \cite{Sinopoli}. This is in contrast to the work \cite{ramesh2012design} where the control signal can affect the states at future times which leads to a dual effect.

\subsection{Erroneous Feedback Communication Channel} \label{ternary_subsec}

In the case of unreliable acknowledgments, the packet loss process $\{\gamma_k, k \geq 0\}$ is not known to the sensor, instead, the sensor receives an imperfect acknowledgment process $\{\hat \gamma_k, k \geq 0\}$ from the receiver. It is assumed that after the transmission of $y_k$ and before transmitting $y_{k+1}$ the sensor has access to the ternary process $\hat \gamma_k\in\{0,1,2\}$ where
\begin{align*}
& \hat \gamma_k=
  \begin{cases}
    0~\textrm{or}~1 &\text{if $\beta_k=1$}\\
    2&\text{if $\beta_k=0$}
  \end{cases}
\end{align*}
with given dropout probability $\eta \in [0,1]$ for the binary process $\{\beta_k: k \geq 0\}$, i.e., $\Prob(\beta_k=0)= \eta \in [0,1]$ for all $k \geq 0$. In case $\beta_k=0$ (i.e., $\hat \gamma_k=2$), no signal is received on the feedback link and this results in an erasure. In case $\beta_k=1$, a transmission error may occur, independent of all other random processes, with probability $\epsilon \in [0,1]$. This transmission error results in the reception of $\hat{\gamma}_k=0$ when $\gamma_k=1$, and $\hat{\gamma}_k=1$ when $\gamma_k=0$. We may write the transition probability matrix of the erroneous feedback channel as a homogeneous Markov process with a $2 \times 3$ transition probability matrix
\begin{align*}
& \mathbb{A} = (a_{ij}) = \left[\begin{array}{ccc}
    (1-\epsilon)(1-\eta) & \epsilon(1-\eta) & \eta  \\
     \epsilon(1-\eta) & (1-\epsilon)(1-\eta) & \eta \\
  \end{array} \right]
\end{align*}
where $a_{ij} : = \Prob(\hat \gamma = j-1 | \gamma= i-1)$ for $i \in \{1,2\}$ and $j \in \{1,2,3\}$. This channel model refers to a generalized erasure channel, namely, a {\em binary erasure channel with errors} 
(see Exercise 7.13 in \cite{cover_thomas}). This model is general in the sense that if we let $\eta=0$ then the ternary acknowledgement process reduces to a binary process with the possibility of only transmission errors, and a standard erasure channel when we set
$\epsilon = 0$. Finally, 
 the case of perfect packet receipt acknowledgments studied in \cite{Nourian_Necsys13} is a special case when $\eta$ and $\epsilon$ above are both set to zero. 

The present situation encompasses, as special cases, situations where no
acknowledgments are available (\emph{UDP-case}) and also cases where
acknowledgments are always available (\emph{TCP-case}), see also\cite{ImerYukselBasar} for a discussion in the context of closed loop control with packet dropouts.

\section{Optimal Transmission Energy Allocation Problems Subject to Energy Harvesting Constraints} \label{optim_sec}

In this section we formulate optimal transmission energy allocation problems in order to minimize the trace of the receiver's expected estimation error covariances (\ref{EH:RiccatiEq}) subject to energy harvesting constraints. Unlike the problem formulation in \cite{Alex_CDC12}, in the model of this paper the optimal energy policies are computed at the sensor which has perfect information about the energy harvesting and instantaneous battery levels but has imperfect state information about the packet receipt acknowledgments. 

We consider the realistic scenario of causal information case where the unpredictable future wireless fading channel gains and energy harvesting information are not a priori known to the transmitter. More precisely, the information available at the sensor at any time $k \geq 1$ is given by
\begin{align*}
& \mathcal{I}_k = \{s_t:=(\hat \gamma_{t-1},g_t,H_t,B_t) : 1 \leq t \leq k\} \cup \mathcal{I}_0
\end{align*}
where $\mathcal{I}_0 := \{g_0,H_0,B_0,P_0\}$ is the initial condition. 

The information $\mathcal{I}_k$ is used at the sensor to decide the amount of transmission energy $u_k$ for the packet loss process. A policy $u_k$ for $k \geq 1$ is feasible if the energy harvesting constraint $0 \leq u_k \leq B_k = \min \{B_{k-1}-u_{k-1}+H_{k},B_{\textrm{max}}\}$ is satisfied. The admissible control set is then given by $\mathcal U := \big\{u_{(\cdot)}: u_k~\textrm{is adapted to sigma-field}~\sigma(\mathcal{I}_k )~\textrm{and}~ 0 \leq u_k \leq B_k~ (a.s.)\big\}.$

The optimization problems are now formulated as Markov decision processes with imperfect state information for the following two cases:
 
(i) {\em Finite-time horizon}:
\begin{equation}
\label{MDP_constrained}
\begin{split}
& \min_{\{u_k: 0 \leq k \leq T-1\}} \sum_{k=0}^{T-1} \mathbb{E} [\textrm{tr}(P_{k+1})] \\
& \textrm{s.t.}~~ 0 \leq u_k \leq B_k ~~ (a.s.) \quad 0 \leq k \leq T-1
\end{split}
\end{equation}

and (ii) {\em Long term average (infinite-time horizon)}:
\begin{equation}
\label{MDP_constrained:IH}
\begin{split}
& \min_{\{u_k: k \geq 0\}} \lim\sup_{T \rightarrow \infty} \frac{1}{T} \sum_{k=0}^{T-1} \mathbb{E} [\textrm{tr}(P_{k+1})] \\
& \textrm{s.t.}~~ 0 \leq u_k \leq B_k ~~ (a.s.) \quad k \geq 0
\end{split}
\end{equation}
where $B_k$ is the stored battery energy available at time $k$ which satisfies the battery dynamics (\ref{EHM}). It is evident that the transmission energy at time $k$, $u_k$, affects the amount of stored energy $B_{k+1}$ available at time $k+1$ which in turn affects the transmission energy $u_{k+1}$ since $0 \leq u_{k+1} \leq B_{k+1} = \min \{B_{k}-u_{k}+H_{k+1}, B_{\textrm{max}}\}$ by (\ref{EHM}). In the special case of perfect packet receipt acknowledgments from receiver to sensor, the reader is referred to \cite{Alex_CDC12} for a similar long term average cost formulation under an average transmission power constraint which is a soft constraint unlike the energy harvesting constraint considered here, which is a hard constraint in an almost sure sense.

We note that the expectations in (\ref{MDP_constrained}) and (\ref{MDP_constrained:IH}) are computed over random variables $\{g_k\}$, $\{H_k\}$ and $\{\hat \gamma_k\}$ for given initial condition $\mathcal{I}_0$. Since these expectations are conditioned on the transmission success process of the feedback channel $\{\hat \gamma_k\}$ instead of the packet loss process of the forward channel $\{\gamma_k  \}$, these formulations fall within the general framework of stochastic control problems with imperfect state information. 

It is known that Kalman filtering with packet losses may have unbounded expected estimation error covariances in certain situations (see \cite{Sinopoli}). We now aim to provide sufficient conditions under which the infinite horizon stochastic control problem (\ref{MDP_constrained:IH}) is well-posed in the sense that an exponential boundedness condition for the expected estimation error covariance is satisfied. The reader is referred to \cite{Quevedo_Automatica} for the problem of determining the minimum average energy required for guaranteeing the stability of the Kalman filtering with the packet reception probabilities (\ref{h_defn}) subject to an average sum energy constraint.

Let $\mathbb{G}$ and $\mathbb{H}$ be the time-invariant probability transition laws of the Markovian channel fading process $\{g_k\}$ and the Markovian harvested energy process $\{H_k\}$, respectively. 

We introduce the following assumption:

({\bf A2}) The channel fading process $\{g_k\}$, harvested energy process $\{H_k\}$ and the maximum battery storage $B_{\textrm{max}}$ satisfy the following:
\begin{align}
& \sup_{(g,H)}  \int_{g_k} \! \int_{H_k} \! (1\!-\!h(g_k \min\{H_k, B_{max}\})) \Prob(g_{k}|g_{k-1}\!=\!g) \notag \\
& \hspace{1cm} \times \Prob(H_k| H_{k-1}\!=\!H)  dg_k dH_k  \leq \frac{\rho}{||A||^2}, \quad k \geq 0 \label{A2-Exp}
\end{align}
for some $\rho \in [0,1)$. \qed

\begin{thm}\label{Stability:thm} Assume ({\bf A2}) holds. Then there exist energy allocations $\{u_k\}$ such that $\{P_k\}$ in (\ref{EH:RiccatiEq}) is exponentially bounded in norm, i.e., 
\begin{align} \label{ExpSta}
& \E \|P_k\| \leq \alpha \rho^k + \beta, \quad k \geq 0
\end{align} 
for some non-negative scalars $\alpha$ and $\beta$. As a result, the stochastic optimal control problem (\ref{MDP_constrained:IH}) is well-posed.
\end{thm}

\textit{Proof}: Based on Theorem 1 in \cite{Quevedo_TAC13}, a sufficient condition for exponential stability in the sense of (\ref{ExpSta}) is that 
\begin{align*}
&\sup_{(g,H)} \Prob(\gamma_k=0 | g_{k-1}=g, H_{k-1}=H) 
\\& = \!\! \sup_{(g,H)} \int_{g_k} \! \int_{H_k} \! \Prob(\gamma_k \!=\!0 | g_k \!=\! g' \!, \! H_{k}\!=\!H' \!, \!g_{k-1}\!=\!g, \! H_{k-1}\!=\!H) \\
& \hspace{2cm} \times \Prob(g_{k}, H_{k}|g_{k-1}\!=\!g, H_{k-1}\!=\!H)  dg_k dH_k  \\
& =\!\! \sup_{(g,H)} \int_{g_k} \! \int_{H_k} \! \Prob(\gamma_k \!=\!0 | g_k \!=\! g' \!, \! H_{k}\!=\!H' \!, \!g_{k-1}\!=\!g, \! H_{k-1}\!=\!H) \\
&\hspace{2cm} \times \Prob(g_{k}|g_{k-1}\!=\!g) \Prob(H_k| H_{k-1}\!=\!H)  dg_k dH_k\\
& = \sup_{(g,H)} \int_{g_k} \int_{H_k}  (1-h(g_k u_k))\Prob(g_{k}|g_{k-1}=g) \\
& \hspace{2cm} \times  \Prob(H_k| H_{k-1}=H)  dg_k dH_k \leq \frac{\rho}{\|A\|^2} 
\end{align*}
 for some $ \rho \in [0,1)$. We now consider a suboptimal solution scheme to the stochastic optimal control problem (\ref{MDP_constrained:IH}) where the full amount of energy harvested at each time step is used, i.e., $u_0=B_0$ and $u_k = \min\{H_k, B_{\textrm{max}}\}$ for $k \geq 1$. Then (\ref{A2-Exp})
will be a sufficient condition in terms of the channel fading process, harvested energy process and the maximum battery storage. Therefore, Assumption ({\bf A2}) provides a sufficient condition for the exponential boundedness (\ref{ExpSta}) of the expected estimation error covariance. \qed

\begin{remark}
The condition (\ref{A2-Exp}) given by Assumption ({\bf A2}) may not be easy to verify for all values of $g$, $H$ and $k$. If we assume that the channel fading and harvested energy processes are stationary then it won't be necessary to verify the condition for all $k$. Furthermore, in the two most commonly used models of i.i.d. processes and finite state Markov chains, the condition can be simplified as follows: 

(i) If $\{g_k\}$ and $\{H_k\}$ are i.i.d., (\ref{A2-Exp}) yields
\begin{align*}
& \int_{g_k} \! \int_{H_k}  (1-h(g_k \min\{(H_k, B_{max}\})) \Prob(g_{k}) \Prob(H_k)  dg_k dH_k  \\
& \hspace{3cm} \leq \frac{\rho}{||A||^2}.
\end{align*}

(ii) If $\{g_k\}$ and $\{H_k\}$ are stationary finite state Markov chains with $M$ and $N$ states respectively, (\ref{A2-Exp}) yields
\begin{align*}
& \max_{(i,j)} \sum_{i'=1}^M \sum_{j'=1}^N (1-h(i \min\{j,B_{\textrm{max}}\})) \Prob(g_k=i'|g_{k-1}=i) \\
& \hspace{3cm} \times \Prob(H_k=j'|H_{k-1}=j)  \leq \frac{\rho}{\|A\|^2}.
\end{align*} \qed
\end{remark} 

\section{Solutions to the Optimal Transmission Energy Allocation Problems Via Dynamic Programming} \label{optimSol_sec}

The stochastic control problems (\ref{MDP_constrained}) and (\ref{MDP_constrained:IH}) can be regarded as Markov Decision Process (MDP) \cite{altman1999constrained} problems with imperfect state information \cite{kumar1986stochastic,bertsekas1995dynamic}. In these formulations the energy harvesting sensor does not have perfect knowledge about whether its transmissions have been received at the receiver or not due to the existence of an imperfect feedback communication channel. Hence, at time $k$ the sensor has only ``imperfect state information'' about $\{P_t: 1 \leq t \leq k\}$ via the acknowledgment process $\{\hat \gamma_t, 0 \leq t \leq k-1 \}$. In this section we reduce the stochastic control problems with imperfect state information (\ref{MDP_constrained}) and (\ref{MDP_constrained:IH}) to ones with perfect state information by using the notion of information-state \cite{kumar1986stochastic}. 

\subsection{Information-State Dynamics}

For $k \geq 0$ denote 
\begin{align*}
& z^k := \{P_0, \hat \gamma_0, \dots, \hat \gamma_{k},g_0,\dots,g_{k-1},\dots,u_0,\dots,u_{k-1}\}
\end{align*}
as all observations about the receiver's Kalman filtering state estimation error covariance at the sensor after the transmission of $y_k$ and before transmitting $y_{k+1}$. We set $z^{-1}:=\{P_0\}$. The so-called {\it information-state} is defined by
\begin{align}
& f_{k+1}(P_{k+1}|z^k,g_k,u_k) = \Prob(P_{k+1}|z^k,g_k,u_k), \quad k \geq 0
\end{align}
which is the conditional probability of estimation error covariance $P_{k+1}$ given $z^k$, $g_k$ and $u_k$. The following lemma shows how $f_{k+1}(\cdot|z^k,g_k,u_k)$ can be determined from $f_k(\cdot|z^{k-1},g_{k-1},u_{k-1})$ together with $\hat \gamma_k$, $g_k$ and $u_k$.

\begin{lem} \label{lem:information-state} The information-state $f_{(\cdot)}$ satisfies the following dynamics
\begin{align} 
& f_{k+1}(P_{k+1}|z^k,g_k,u_k) =  \sum_{\gamma_k \in\{0,1\}}  \!\!\! \Big[ \int_{P_k} \Big(\Prob(P_{k+1} | P_k,\gamma_k) \notag \\
& \hspace{0.5cm} \times f_{k} (P_k|z^{k-1},g_{k-1},u_{k-1}) \Big)dP_k \notag \\
& \hspace{0.5cm}   \times \frac{\Prob(\hat \gamma_k|\gamma_k) \times \Prob(\gamma_k | g_k,u_k)}{ \sum_{\gamma_k \in \{0,1\}}\Prob(\hat \gamma_k|\gamma_k) \times \Prob(\gamma_k|g_k,u_k)} \Big], \quad k \geq 0 \label{InfState_recur_tot}
\end{align}
with $f_0(P_0|z^{-1}) = \delta(P_0)$ where $\delta$ is the Dirac delta function. 
\end{lem}

{\it Proof}: See the Appendix. \qed

It is important to note that the information-state dynamics (\ref{InfState_recur_tot}) depends on the fading channel gains $\{g_k\}$ and sensor transmission energy allocation policies $\{u_k\}$ via the packet reception probabilities (\ref{h_defn}). Hence, we may write (\ref{InfState_recur_tot}) as
\begin{align} 
& f_{k+1}(P_{k+1}|z^k,g_k,u_k) \notag \\
& \hspace{0.5cm} = \Phi_k\big[f_{k}(\cdot|z^{k-1},g_{k-1},u_{k-1}),\hat \gamma_k,g_k,u_k\big](P_{k+1}) \label{InfState_recur_tot:2}
\end{align}
for $k \geq 0$. Note that $\Phi_k$ in (\ref{InfState_recur_tot:2}) depends on the entire function $f_{k}(\cdot|z^{k-1})$ and not just its value at any particular $P_{k}$. 

In the following sections the stochastic control problems with imperfect state information (\ref{MDP_constrained}) and (\ref{MDP_constrained:IH}) are reduced to problems with perfect state information where the state is given by the information-state $f_{(\cdot)}$. The resulting stochastic problems with perfect information are approached via the dynamic programming principle. 

We establish some notation. Let the binary random variable $\gamma$ be defined akin to $\gamma_k$ in (\ref{EH:RiccatiEq}), then for a given $P$ denote
\begin{align} \label{EH:RiccatiEqOperator}
& \! \! \!  \mathcal L (P,\gamma) \!  := \! A P A^T \!\!+\! Q \!-\! \gamma A P C^T\! (C P C^T \!+\! R)^{-1} C P A^T
\end{align}
as the random Riccati equation operator. Let $\mathcal S_{+}^n$ be the set of all $(n \times n)$ nonnegative definite matrices. Then, we denote the space of all probability density functions on $\mathcal S_{+}^n$ as $\Pi$ where $\int_{\mathcal S_{+}^n} \pi(P) dP =1$ for any $\pi \in \Pi$. Let the ternary random variable $\hat \gamma$ be defined akin to $\hat \gamma_k$ in Section \ref{ternary_subsec}. Then, based on the information-state recursion (\ref{InfState_recur_tot:2}) denote 
\begin{align} 
& \tilde \pi = \Phi \big[\pi,\hat \gamma,g,u\big] := \!\!\!\!\!\! \sum_{\gamma \in\{0,1\}} \!\!\! \Big[ \int_{P} \Prob \big(\mathcal L(P,\gamma) \big| P,\gamma\big) \pi(P) dP \notag \\
& \hspace{1cm} \times \frac{\Prob(\hat \gamma|\gamma) \times \Prob(\gamma|g,u)}{ \sum_{\gamma \in \{0,1\}}\Prob(\hat \gamma|\gamma) \times \Prob(\gamma|g,u)} \Big] \label{FunPhi}
\end{align}
for given $\pi \in \Pi$, fading channel gain $g$ and sensor transmission energy allocation $u$. 

\begin{remark} \label{perffeedRemark} In the special case of perfect packet receipt acknowledgments, where $\eta$ and $\epsilon$ in Section \ref{ternary_subsec} are set to zero, the problems (\ref{MDP_constrained}) and (\ref{MDP_constrained:IH}) become stochastic control problems with perfect state information. In this case the probability density functions $\pi$ and $\tilde \pi$ in the information-state recursion (\ref{FunPhi}) become Dirac delta functions. \qed
\end{remark}

\subsection{Dynamic Programming Principle}  \label{Sec:CIC}

In this section, the transmission energy allocation policy is computed offline from the Bellman dynamic programming equations given below. 

Some notation is now presented. Given the fading channel gain $g$ and the harvested energy $H$ at time $k \geq 0$ we denote the corresponding fading channel gain and the harvested energy at time $k+1$ by $\tilde g$ and $\tilde H$, respectively. We recall that both fading channel gains $\{g_k\}$ and harvested energies $\{H_k\}$ are modelled as first-order homogeneous Markov processes (see Section \ref{model_sec}).   

\subsubsection{Finite-Time Horizon Bellman Equation} The imperfect state information stochastic control problem (\ref{MDP_constrained}) is solved in the following Theorem.
\begin{thm} \label{Causal:FH:Bellman} For given initial condition $\mathcal{I}_0 = \{g_0,H_0,B_0,P_0\}$ the value of the finite-time horizon minimization problem (\ref{MDP_constrained}) is given by $V_0(\mathcal{I}_0)$ which can be computed recursively from the backward Bellman dynamic programming equation
\begin{align}
& V_k (\pi,g,H,B) = \min_{0 \leq u \leq B} \Big\{\mathbb{E} \big[\tr\big(\mathcal L(P,\gamma)\big)\big|\pi,g,u\big] \notag \\
& \hspace{0.5cm} + \mathbb{E} \Big[V_{k+1} \big(\Phi \big[\pi,\hat \gamma,g,u\big],\tilde g,\tilde H,\min\{B-u+\tilde H, B_{\textrm{max}}\}\big) \notag \\
& \hspace{3cm} \big| \pi,g,H,u\Big] \Big\} \label{perfecfeed_dp}, \qquad 0 \leq k \leq T-1
\end{align}
where $\pi \in \Pi$. The terminal condition is given as 
\begin{align*}
& V_{T}(\pi,g,H,B):= \min_{0 \leq u \leq B} \mathbb{E} \big[\tr\big(\mathcal L(P,\gamma)\big)\big|\pi,g,u\big] \\
& \hspace{2.2cm} = \mathbb{E} \big[\tr\big(\mathcal L(P,\gamma)\big)\big|\pi,g,B\big]
\end{align*} 
where all available energy is used for transmission in the final time $T$. 
\end{thm}

\textit{Proof}: The proof follows from the dynamic programming principle for stochastic control problems with imperfect state information (see Theorem 7.1 in \cite{kumar1986stochastic}). \qed

Based on Remark \ref{perffeedRemark}, it is important to note that in the special case of perfect packet receipt acknowledgments, where $\eta$ and $\epsilon$ in Section \ref{ternary_subsec} are set to zero, the Bellman equation (\ref{perfecfeed_dp}) is written with respect to Dirac delta functions in space $\Pi$, i.e., $\pi(\cdot) =\delta_{(\cdot)} \in \Pi$ (see Section 4 in \cite{Nourian_Necsys13}).

The solution to the imperfect state information stochastic control problem (\ref{MDP_constrained}) is then given by
\begin{align}
& u_k^o(\pi,g,H,B) = \arg \min_{0 \leq u \leq B} \Big\{\mathbb{E} \big[\tr\big(\mathcal L(P,\gamma)\big)\big|\pi,g,u\big] \notag \\
& \hspace{0.5cm} + \mathbb{E} \Big[V_{k+1} \big(\Phi \big[\pi,\hat \gamma,g,u\big],\tilde g,\tilde H,\min\{B-u+\tilde H, B_{\textrm{max}}\}\big) \notag \\
& \hspace{3cm} \big| \pi,g,H,u\Big] \Big\}  \label{OptimalPoli}, \qquad 0 \leq k \leq T-1
\end{align}
with $ u_T^o(\pi,g,H,B) = B$, where $V_{k+1}(\cdot)$ is the solution to the Bellman equation (\ref{perfecfeed_dp}).

For computational purposes, we now simplify the terms in (\ref{perfecfeed_dp}). First, we have
\begin{align*}
& \mathbb{E} \big[\tr\big(\mathcal L(P,\gamma)\big)\big|\pi,g,u\big]  = \int_{P} \tr \big(A P A^T + Q\big) \pi (P) dP \notag\\
& \hspace{0.1cm} - h(gu) \times \!\! \int_{P} \!\! \tr \Big(A P C^T [C PC^T + R]^{-1} C P A^T \Big) \pi(P) dP
\end{align*} 
with the constraint that $0 \leq u \leq B$. Since the mutually independent processes $\{g_{k}\}$ and $\{H_{k}\}$ are independent of other processes and random variables, we may write
\begin{align}
& \mathbb{E} \Big[V_{k+1} \big(\Phi \big[\pi,\hat \gamma,g,u\big],\tilde g,\tilde H,\tilde B\big) \big| \pi,g,H,u\Big] \Big\} \
\notag \\ 
& \hspace{0.2cm} =  \int\limits_{\tilde g, \tilde H} V_{k+1} \Big(\Phi [\pi,\hat \gamma,g,u],\tilde g,\tilde H,\tilde B\Big) \notag \\
& \hspace{4cm} \times \mathbb{G}(\tilde g|g) ~ \mathbb{H} (\tilde H|H) d \tilde g d \tilde H
\label{ValueMarkovExp}
\end{align}
where $\tilde B:= \min\{B-u+\tilde H, B_{\textrm{max}}\}$, and $\mathbb{G}$ and $\mathbb{H}$ are the probability transition laws of the Markovian processes $\{g_k\}$ and $\{H_k\}$, respectively. But, 
\begin{align*}
& E\big[\Phi [\pi,\hat \gamma,g,u] \big|\pi,g,u\big] = \Prob(\hat \gamma=0) \times \Phi [\pi,0,g,u] \\
& \hspace{0.4cm} +  \Prob(\hat \gamma=1) \times \Phi [\pi,1,g,u] + \Prob(\hat \gamma=2) \times \Phi [\pi,2,g,u]
\end{align*}
where the function $\Phi$ is defined in (\ref{FunPhi}).

\begin{remark}
The expression (\ref{ValueMarkovExp}) can be simplified further in the two following cases: \\
(i) If $\{g_k\}$ and $\{H_k\}$ are i.i.d., then the right hand side term in (\ref{ValueMarkovExp}) becomes
\begin{align*}
& \int\limits_{\tilde g, \tilde H} V_{k+1} \Big(\Phi [\pi,\hat \gamma,g,u],\tilde g,\tilde H,\tilde B \Big) \Prob (\tilde g) \Prob(\tilde H) d \tilde g d \tilde H
\end{align*}
where $\tilde B= \min\{B-u+\tilde H, B_{\textrm{max}}\}$. \\
(ii) If $\{g_k\}$ and $\{H_k\}$ are finite state Markov chains with $M$ and $N$ states respectively, then the right term in (\ref{ValueMarkovExp}) becomes
\begin{align*}
&  \sum_{i=1}^M \sum_{j=1}^N V_{k+1} \Big(\Phi [\pi,\hat \gamma,g,u], i,j,\tilde B(j)\Big) \times (\Prob(g) \mathbb{G})_i (\Prob(H) \mathbb{H})_j
\end{align*}
where $\tilde B(j) := \min\{B-u+j, B_{\textrm{max}}\}$, $\Prob(g) :=[\Prob(g=1) ~  \dots ~ \Prob(g=M)]$, $\Prob(H) :=[\Prob(H=1)~ \dots ~\Prob(H=N)]$, $\mathbb{G}$ and $\mathbb{H}$ are the probability transition matrices for $\{g_k\}$ and $\{H_k\}$, respectively, and $(\Prob(g) \mathbb{G})_i $ denotes the i-th component of the vector $\Prob(g) \mathbb{G}$. \qed
\end{remark}

Note that the solution to the dynamic programming equation can only be obtained numerically and there is no closed form solution. In fact, even for a horizon 2 problem with causal information and perfect feedback acknowledgment, it can be shown that the optimal solution cannot be obtained in closed form. It can be observed however that for a fixed battery level, the energy allocation generally increases with the channel gain and when the channel gain is above some threshold, all of the available battery energy is used for transmission.
Similarly, when the channel gain is kept fixed, the energy allocation is equal to the available energy and increases with increasing battery energy level. Although after some point, the energy allocated for transmission becomes less than the available energy and some energy is saved for future transmissions.
\subsubsection{Long Term Average (Infinite-Time Horizon) Bellman Equation} We present the solution to the imperfect state information stochastic control problem (\ref{MDP_constrained:IH}) in the following Theorem.
\begin{thm} \label{Causal:IH:Bellman} Independent of the initial condition $\mathcal{I}_0 = \{g_0,H_0,B_0,P_0\}$, the value of the infinite-time horizon minimization problem (\ref{MDP_constrained:IH}) is given by $\rho$ which is the solution of the average-cost optimality (Bellman) equation
\begin{align}
& \rho + V(\pi,g,H,B) = \min_{0 \leq u \leq B} \Big\{\mathbb{E} \big[\tr\big(\mathcal L(P,\gamma)\big)\big|\pi,g,u\big] \notag \\
& \hspace{0.5cm} + \mathbb{E} \Big[V\big(\Phi \big[\pi,\hat \gamma,g,u\big],\tilde g,\tilde H,\min\{B-u+\tilde H, B_{\textrm{max}}\}\big) \notag \\
&\hspace{3cm} \big| \pi,g,H,u\Big] \Big\} \label{perfecfeed_dp:IH},
\end{align}
where $\pi \in \Pi$, and $V$ is called the relative value function.
\end{thm}

\textit{Proof}: See the Appendix. \qed

The stationary solution to the imperfect state information stochastic control problem (\ref{MDP_constrained:IH}) is then given by
\begin{align}
& u^o(\pi,g,H,B) = \arg \min_{0 \leq u \leq B} \Big\{\mathbb{E} \big[\tr\big(\mathcal L(P,\gamma)\big)\big|\pi,g,u\big] \notag \\
& \hspace{0.5cm} + \mathbb{E} \Big[V\big(\Phi \big[\pi,\hat \gamma,g,u\big],\tilde g,\tilde H,\min\{B-u+\tilde H, B_{\textrm{max}}\}\big) \notag \\
& \hspace{3cm} \big| \pi,g,H,u\Big] \Big\} \label{Optimal:Sol}
\end{align}
where $V(\cdot)$ is the solution to the average cost Bellman equation (\ref{perfecfeed_dp:IH}).

\begin{remark} Equation (\ref{perfecfeed_dp:IH}) together with the control policy $u^o$ defined in (\ref{Optimal:Sol}) is known as the average cost optimality equations. If a control $u^o$, a measurable function $V$, and a constant $\rho$ exist which solve equations (\ref{perfecfeed_dp:IH})-(\ref{Optimal:Sol}), then the control $u^o$ is optimal, and $\rho$ is the optimal cost in the sense that
\begin{equation*}
\begin{split}
& \lim\sup_{T \rightarrow \infty} \frac{1}{T} \sum_{k=0}^{T-1} \mathbb{E} [\textrm{tr}(P_{k+1})|u^o] = \rho
\end{split}
\end{equation*}
and for any other control policy $\{u_k: k \geq 0\}$ such that $0 \leq u_k \leq B_k$, a.s., we have \begin{equation*}
\begin{split}
& \lim\sup_{T \rightarrow \infty} \frac{1}{T} \sum_{k=0}^{T-1} \mathbb{E} [\textrm{tr}(P_{k+1})|u] \geq \rho.
\end{split}
\end{equation*}
The reader is referred to\cite{arapostathis1993discrete} for a proof of the average cost optimality equations and related results. \qed
\end{remark} 

We note that discretized versions of the Bellman equations (\ref{perfecfeed_dp}) or (\ref{perfecfeed_dp:IH}), which in particular includes the discretization of the space of probability density functions $\Pi$, is used for the numerical computation to find suboptimal solutions to the stochastic control problems (\ref{MDP_constrained}) and (\ref{MDP_constrained:IH}). As the number of discretization levels increases, it is expected that these discretized (suboptimal) solutions converge to the optimal solutions \cite{YuBertsekas}. We solve the Bellman equations (\ref{perfecfeed_dp}) and (\ref{perfecfeed_dp:IH}) by the use of {\it value iteration} and {\it relative value iteration} algorithms, respectively (see Chapter 7 in \cite{bertsekas1995dynamic}). 

\begin{remark}
The causal information pattern is clearly relevant to the most practical scenario. However, it is also instructive to consider the non-causal information scenario where the sensor has a priori information about the energy harvesting $\{H_k\}$ process and the fading channel gains $\{g_k\}$ for all time periods including the future ones. This may be feasible in the situation of known environment where the wireless channel fading gains and the harvested energies are predictable \cite{ho2012optimal}. More importantly, the performance of the non-causal information case can serve as a benchmark (a lower bound) for the causal case. Indeed, we present some performance comparison between the performances in the causal and the non-causal case in the {\em Numerical Examples} section. Note that the energy allocation problems for the non-causal case can be solved using similar techniques to Section \ref{Sec:CIC}, and the details are omitted for brevity. \qed
\end{remark}

\section{Suboptimal Transmission Energy Allocation Problems and Their Solutions} \label{Sub_sec}

The optimal solutions presented in Section \ref{optim_sec} require us to compute the solution of Bellman equations in the space of probability densities $\Pi$. In this section we consider the design of suboptimal policies which are computationally much less intensive than the optimal solutions of Section \ref{optimSol_sec}. 

Here, we only present suboptimal solutions to the finite-time horizon stochastic control problem (\ref{MDP_constrained}). Following the same arguments one can design similar suboptimal solutions to the infinite-time horizon problem.

In this case we formulate the problem of minimizing the expected estimation error covariance as
\begin{equation}
\label{MDP_constrained:SubOp}
\begin{split}
& \min_{\{u_k: 0 \leq k \leq T\}} \sum_{k=0}^{T-1} \mathbb{E} \big[\textrm{tr}\big(\hat P_{k+1}) \big| \{\hat \gamma_l\}_{l=0}^{k-1},\{u_l\}_{l=0}^{k},P_0\big] \\
& \hspace{0.5cm}  \equiv \min_{\{u_k: 0 \leq k \leq T\}} \sum_{k=0}^{T-1} \mathbb{E} \big[\textrm{tr}\big(\hat P_{k+1}) \big| \hat P_k, u_k\big]  \\
& \textrm{s.t.}~~ 0 \leq u_k \leq B_k ~~ (a.s.) \quad 0 \leq k \leq T-1
\end{split}
\end{equation}
where $\hat P_{(\cdot)}$ is an estimate of $P_{(\cdot)}$ computed by the sensor based on the following recursive equations (with $\hat P_0  = P_0$): 

(i) In the case $\hat \gamma_k = 0$ we have
\begin{align*}
& \hat P_{k+1} \!:=\!\! \big(A \hat P_k A^T + Q \big) \!\! \times\!\! \frac{\Prob(\hat \gamma_k=0|\gamma_k=0) \times \Prob(\gamma_k=0)}{ \sum_{\gamma_k \in \{0,1\}}\Prob(\hat \gamma_k=0|\gamma_k) \times \Prob(\gamma_k)} \\
& \hspace{0.5cm} + \big(A \hat P_k A^T + Q - A \hat P_k C^T [C\hat P_k C^T+ R]^{-1} C \hat P_k A^T \big) \\
& \hspace{1.5cm} \times \frac{\Prob(\hat \gamma_k=0|\gamma_k=1) \times \Prob(\gamma_k=1)}{ \sum_{\gamma_k \in \{0,1\}}\Prob(\hat \gamma_k=0|\gamma_k) \times \Prob(\gamma_k)}. 
\end{align*}

(ii) in the case $\hat \gamma_k = 1$ we have
\begin{align*}
& \hat P_{k+1} \!:=\!\! \big(A \hat P_k A^T + Q \big)  \!\! \times\!\!  \frac{\Prob(\hat \gamma_k=1|\gamma_k=0) \times \Prob(\gamma_k=0)}{ \sum_{\gamma_k \in \{0,1\}}\Prob(\hat \gamma_k=1|\gamma_k) \times \Prob(\gamma_k)} \\
& \hspace{0.5cm} + \big(A \hat P_k A^T + Q - A \hat P_k C^T [C\hat P_k C^T+ R]^{-1} C \hat P_k A^T \big) \\
& \hspace{1.5cm} \times \frac{\Prob(\hat \gamma_k=1|\gamma_k=1) \times \Prob(\gamma_k=1)}{ \sum_{\gamma_k \in \{0,1\}}\Prob(\hat \gamma_k=1|\gamma_k) \times \Prob(\gamma_k)}. 
\end{align*}

(iii) In the case $\hat \gamma_k = 2$ we have
\begin{align*}
& \hat P_{k+1} := A \hat P_k A^T + Q  - \Prob(\gamma_k =1) \\
& \hspace{3cm} \times A \hat P_k C^T  [C\hat P_k C^T+ R]^{-1} C \hat P_k A^T.
\end{align*}

The reason that the solution to the stochastic control problem (\ref{MDP_constrained:SubOp}) is called suboptimal is that the true estimation error covariance matrix $P_{(\cdot)}$ in (\ref{EH:RiccatiEq}) is replaced by its estimate $\hat P_{(\cdot)}$. The intuition behind these recursive equations can be explained as follows.
Note that in the case of perfect feedback acknowledgements, the error covariance is updated as $P_{k+1} = A P_k A^T + Q$ in case $\gamma_k=0$, and $P_{k+1} = A P_k A^T + Q - A P_k C^T (C P_k C^T + R)^{-1} C P_k A^T $ in case $\gamma_k=1$. In our imperfect acknowledgement model, even when it is received, errors can occur such that $\hat{\gamma}_k=0$ is received when $\gamma_k=1$, and $\hat{\gamma}_k=1$ is received when $\gamma_k=0$. Thus the recursions given in (i) and (ii) are the weighted (by the corresponding error event probabilities) combinations of the error covariance recursions in the case of perfect feedback acknowledgements. In the case $\hat{\gamma}_k=2$ where an erasure 
occurs, taking the average of the error covariances in the cases $\gamma_k=0$ and $\gamma_k=1$ is intuitively a reasonable thing to do, which motivates the recursion in (iii).

Note that $\Prob(\hat \gamma_k)= \sum_{\gamma_k \in \{0,1\}} \Prob(\hat \gamma_k | \gamma_k) \Prob(\gamma_k)$ where the conditional probabilities are given in Section \ref{ternary_subsec}. This together with the recursive equations of $\hat P_{(\cdot)}$ yields
\begin{align}
& \mathbb{E} [\hat P_{k+1}|\hat P_k,g_k,u_k] = A \hat P_k A^T + Q - h(g_ku_k) \notag\\
& \hspace{0.1cm} \times \big(A \hat P_k A^T + Q -  A \hat P_k C^T  [C\hat P_k C^T+ R]^{-1} C \hat P_k A^T \big)
\end{align}

Since the expression $\mathbb{E} [\hat P_{k+1}| \hat P_k,g_k,u_k]$ is of the same form as $\mathbb{E} [P_{k+1}| P_k, g_k,u_k] $ when $P_k$ is replaced by $\hat P_k$, the Bellman equation for problem (\ref{MDP_constrained:SubOp}) is given by a similar equation to the case of perfect feedback communication channel considered in \cite{Nourian_Necsys13} which is presented in the following theorem.

\begin{thm} \label{Causal:FH:Bellman:SubOp} For given initial condition $\mathcal{I}_0 = \{g_0,H_0,B_0,P_0\}$ the value of the finite-time horizon minimization problem (\ref{MDP_constrained:SubOp}) is given by $V_0(\mathcal{I}_0)$ which can be computed recursively from the backward Bellman dynamic programming equation
\begin{align}
& V_k (\hat{P},g,H,B) = \min_{0 \leq u \leq B} \Big\{\mathbb{E} \big[\tr\big(\mathcal L(\hat{P},\gamma)\big)\big|\hat{P},g,u\big] \notag \\
& \hspace{0.5cm} + \mathbb{E} \Big[V_{k+1} \big(\mathcal L(\hat{P},\gamma),\tilde g,\tilde H,\min\{B-u+\tilde H, B_{\textrm{max}}\}\big) \notag \\
& \hspace{2.8cm} \big|\hat{P},g,H,u\Big] \Big\} \label{perfecfeed_dp:SubOpt}, \qquad 0 \leq k \leq T-1
\end{align}
with terminal condition
\begin{align*}
& V_{T}(\hat{P},g,H,B):= \min_{0 \leq u \leq B} \mathbb{E} \big[\tr\big(\mathcal L(\hat{P},\gamma)\big)\big|\hat{P},g,u\big]  \\
& \hspace{2.2cm} = \mathbb{E} \big[\tr\big(\mathcal L(\hat{P},\gamma)\big)\big|\hat{P},g,B\big]
\end{align*} 
where all available energy is used for transmission in the final time $T$. \qed
\end{thm}

The solution to the stochastic control problem (\ref{MDP_constrained:SubOp}) which is a suboptimal solution to (\ref{MDP_constrained}) is given by
\begin{align*}
& u_k^o(\hat{P},g,H,B) = \arg \min_{0 \leq u \leq B} \Big\{\mathbb{E} \big[\tr\big(\mathcal L(\hat{P},\gamma)\big)\big|\hat{P},g,u\big] \notag \\
& \hspace{0.4cm} + \mathbb{E} \Big[V_{k+1} \big(\mathcal L(\hat{P},\gamma),g,u\big],\tilde g,\tilde H,\min\{B-u+\tilde H, B_{\textrm{max}}\}\big) \big| \\
& \hspace{2.8cm} \hat{P},g,H,u\Big] \Big\}, \qquad 0 \leq k \leq T-1
\end{align*}
with $u_T^o(\hat{P},g,H,B) = B$, where $V_{k+1}(\cdot)$ is the solution to the Bellman equation (\ref{perfecfeed_dp:SubOpt}).

\section{Some Structural Results on the Optimal Energy Allocation Policies} \label{Str_sec}
In this section the structure of the optimal transmission energy allocation policies (\ref{OptimalPoli}) is studied for the case of the finite-time horizon stochastic control problem (\ref{MDP_constrained}). Following the same arguments one can show similar structural results for the infinite-time horizon problem (\ref{MDP_constrained:IH}).

\begin{lem}\label{convex:lem} Assume $h(\cdot)$ in (\ref{h_defn}) is a concave function in $u$ given $g$. Then, given $\pi, g$ and $H$, the value function $V_k(\pi,g,H,B)$ in (\ref{perfecfeed_dp}) is convex in $B$ for $0 \leq k \leq T$. As a result,
\begin{align*}
& V_0(P_0,g_0,H_0,B_0) = \min_{\{0 \leq u_k \leq B_k\}_{k=0}^{T-1}} \sum_{k=0}^{T-1} \mathbb{E} [\textrm{tr}(P_{k+1})]
\end{align*}
is convex in $B_0$.   
\end{lem}    

\textit{Proof}: We let $s:=(\pi, g,H,B)$. First, note that, for given $\pi, g$ and $H$, the final time value function
\begin{align*}
& V_{T}(s) = \min_{0 \leq u \leq B} \mathbb{E} \big[\tr\big(\mathcal L (P,\gamma)\big)\big|\pi,g,u\big] \\
& \hspace{0.9cm} = \mathbb{E} \big[\tr\big(\mathcal L(P,\gamma)\big)\big|\pi,g,B\big]
\end{align*}
is a convex function in $B$ due to the fact that $h(\cdot)$ is a concave function in $u$ given $g$ (see Lemma 2 in \cite{Alex_CDC12}). Now assume that $V_{k+1}(s)$ is convex in $B$ for given $\pi,g$ and $H$. Then, for given $H$ and $u$, the function
\begin{align*}
& V_{k+1}(\pi,g,H,\min \{B-u+H,B_{\textrm{max}}\})
\end{align*}
is convex in $B$, since it is the minimum of $V_{k+1}(\pi,g,H,B_{\textrm{max}})$ which is a constant independent of $B$, and by the induction hypothesis the convex function $V_{k+1}(\pi,g,H,B-u+H)$ in $B$. Since the expectation operator preserves convexity, 
\begin{align*}
& E \Big[V_{k+1} \big(\Phi \big[\pi,\hat \gamma,g,u\big],\tilde g,\tilde H,\min\{B-u+\tilde H, B_{\textrm{max}}\}\big) \\
& \hspace{3cm} \big| \pi,g,H,u\Big] 
\end{align*}
given in (\ref{ValueMarkovExp}) is a convex function in $B$. But, $V_{k}(s)$ in (\ref{perfecfeed_dp}) is the infimal convolution of two convex functions in $B$ for given $\pi,g$ and $H$ and hence is convex in $B$ (see the proof of Theorem 1 in \cite{ho2012optimal}). \qed

We now present the main Theorem of this section which gives structural results on the optimal energy allocation policies (\ref{OptimalPoli}).

\begin{thm}\label{Structure:thm} Assume $h(\cdot)$ in (\ref{h_defn}) is a concave function in $u$ given $g$. Then, given $\pi, g$ and $H$, the optimal transmission energy allocation policy $u^o_k(\pi,g,H,B)$ given in (\ref{OptimalPoli}) is non-decreasing in $B$ for $0 \leq k \leq T$.
\end{thm}

\textit{Proof}: Assume $\pi, g$ and $H$ are fixed. We define
\begin{align*}
& L(B,u) = \mathbb{E} \big[\tr\big(\mathcal L(P,\gamma)\big)\big|\pi,g,u\big] \notag \\
& \hspace{0.3cm} + \mathbb{E} \Big[V_{k+1} \big(\Phi \big[\pi,\hat \gamma,g,u\big],\tilde g,\tilde H,\min\{B-u+\tilde H, B_{\textrm{max}}\}\big) \\
& \hspace{3cm} \big| \pi,g,H,u\Big] 
\end{align*}
from (\ref{perfecfeed_dp}). We aim to show that $L(B,u)$ is submodular in $(B,u)$, i.e., for every $u' \geq u$ and $B' \geq B$,
\begin{align} \label{sub}
& L(B', u') - L(B,u') \leq L(B', u) - L(B,u).
\end{align}
It is evident that $\mathbb{E} \big[\tr\big(\mathcal L(P,\gamma)\big)\big|\pi,g,u\big]$ is submodular in $(B,u)$ since it is independent of $B$. Denote
\begin{align*}
& Z(x) := \mathbb{E} \Big[V_{k+1} \big(\Phi \big[\pi,\hat \gamma,g,u\big],\tilde g,\tilde H,\min\{x+\tilde H, B_{\textrm{max}}\}\big) \\
& \hspace{3cm} \big| \pi,g,H,u\Big]. 
\end{align*}
Since $Z(x)$ is convex in $x$ (by Lemma \ref{convex:lem}) we have 
\begin{align*}
& Z(x+\epsilon)-Z(x) \leq Z(y+\epsilon)-Z(y), \quad x \leq y, ~\epsilon \geq 0
\end{align*}
(see Proposition 2.2.6 in \cite{simchi2004logic}). Now let $x=B-u'$, $y=B-u$ and $\epsilon=B'-B$. Then, we have the submodularity condition (\ref{sub}) for $Z(B-u)$ \cite{ho2012optimal}. Therefore, $L(B,u)$ is submodular in $(B,u)$. Note that submodularity is a sufficient condition for optimality of monotone increasing policies, i.e., since $L(B, u)$ is submodular in $(B,u)$ then $u^o(B) =\arg\min_u L(B,u)$ is non-decreasing in $B$ (see \cite{topkis2001supermodularity}). \qed 

For fixed $\pi$, $g$ and $H$, let $u^*_k$ be the unique solution to the convex unconstrained minimization problem
\begin{align*}
& u_k^*(\pi,g,H,B) = \arg \min_{u } \Big\{\mathbb{E} \big[\tr\big(\mathcal L(P,\gamma)\big)\big|\pi,g,u\big] \notag \\
& \hspace{0.5cm} + \mathbb{E} \Big[V_{k+1} \big(\Phi \big[\pi,\hat \gamma,g,u\big],\tilde g,\tilde H,\min\{B-u+\tilde H, B_{\textrm{max}}\}\big) \\
& \hspace{3cm} \big| \pi,g,H,u\Big] \Big\} 
\end{align*}
which can be easily solved using numerical techniques such as a bisection search. Then, the structural result of Theorem \ref{Structure:thm} implies that the solution to the constrained problem (\ref{OptimalPoli}) where $0 \leq u \leq B$ will be of the form
\begin{align*}
& u_{k}^o (\pi,g,H,B) = \left\{ \begin{array}{cl}
0  , & \quad \textrm{if} ~ u_k^* \leq 0    \\
u_k^*   , & \quad \textrm{if} ~ 0< u_k^* < B    \\
B , & \quad   \textrm{if} ~  u_k^* \geq B.
\end{array} \right. 
\end{align*}
This also helps to reduce the search space by restricting the search to be in one direction for different $B$ (see the discussion in Section III.C of \cite{ho2012optimal}).

\subsection{Threshold Policy for Binary Energy Allocation Levels}
Note that while solving for the optimal energy allocation level in the Bellman equation requires not only discretization of the state space, but also that of the action space. However, the discretization of the action space to a finite number of energy allocation levels is not often an issue as in practice, a sensor transmitter can be programmed to have a finite number of transmission power/energy levels only. In fact, for simplicity of implementation, often a sensor can be equipped with only two power/energy levels for transmission. Thus it is perfectly natural to consider the scenario where the energy allocation space is binary.
In this section therefore we consider the optimal solution of Section \ref{optimSol_sec} with the assumption
 that the transmission energy allocation control $u_k$ belongs to a two element set $\{E_0,E_1\}$ where $E_0 < E_1$. The monotonicity of Theorem \ref{Structure:thm} yields a threshold structure such that the optimal transmission energy allocation policy is of the form
\begin{align} 
& u_{k}^* (\pi,g,H,B) = \left\{ \begin{array}{cl}
E_0  , & \quad \textrm{if} ~ B \leq B^*(\pi,g,H,B)   \\
E_1  , & \quad  \textrm{otherwise}
\end{array} \right.  \label{THRES:Polic}
\end{align}
for $B \geq E_0$, where $B^*(\pi,g,H,B)$ is the battery storage threshold depending on $\pi$, $g$, $H$ and $B$. This threshold structure simplifies the implementation of the optimal energy allocation significantly. However, it requires the knowledge of the optimal 
battery energy threshold $B^*$ above. In general, there is no 
closed form expression for $B^*$, but it can be found via iterative search algorithms. Here we present a gradient estimate based algorithm
based on Algorithm 1 in \cite{krishnamurthy2012sequential} (after \cite{spall2005introduction}) to find the threshold in the case of the infinite-time horizon formulation (\ref{MDP_constrained:IH}) with perfect packet receipt acknowledgments where $\eta$ and $\epsilon$ in Section \ref{ternary_subsec} are set to zero. 
A similar algorithm can be devised for the imperfect feedback case albeit with increased computational complexity.

First, we establish some notation. Let $V^{(k)}$ be the $k$-th iteration of the relative value algorithm for solving Bellman equation (\ref{perfecfeed_dp:IH}) in the case of perfect feedback. Then, for given $B^*$ and fixed $P$, $g$, $H$ and $B$ denote
\begin{align}
& J^{(k)} (B^*) := \mathbb{E} \big[\tr\big(\mathcal L(P,\gamma)\big)\big|P,g,u^*\big] \notag \\
& \hspace{0.4cm} + \mathbb{E} \big[V^{(k-1)}\big(\mathcal L(P,\gamma),\bar g, \bar H,\min\{B-u^*+\bar H, B_{\textrm{max}}\}\big) \notag \\ 
&\hspace{1.3cm} \big|P,g,H,u^*\big] \label{RVI:Threshold}, \qquad k = 1,2, \cdots
\end{align}
where the threshold policy $u^*$ is defined as
\begin{align*} 
& u^* = \left\{ \begin{array}{cl}
E_0 , & \quad \textrm{if} ~ B \leq B^*   \\
E_1  , & \quad  \textrm{otherwise}.
\end{array} \right.  
\end{align*}
For $n \in \mathbb{N}$, $0.5 < \kappa \leq 1$ and $\omega,\varsigma >0$ we denote $\omega_n:=\frac{\omega}{(n+1)^\kappa}$ and $\varsigma_n := \frac{\varsigma}{(n+1)^\kappa}$. The term $J^{(k)}$ in (\ref{RVI:Threshold}) is the right hand side expression of (\ref{perfecfeed_dp:IH}) without the minimization, where the threshold policy $u^*$, depending on the threshold policy $B^*$, is used in the relative value iteration.

{\it Gradient algorithm for computing the threshold}. For fixed $P$, $g$, $H$ and $B$ in the $k$-th iteration of the relative value algorithm the following steps are carried out:
 
Step 1) Choose the initial battery storage threshold $B^{(0)}$.

Step 2) For iterations $n=0,1, \cdots$
\begin{itemize}
\item Compute the gradient:
\begin{align} 
& \!\! \partial_{B} J^{(k)}_n := \frac{J^{(k)}(B^{(n)}+\omega_n)- J^{(k)}(B^{(n)}-\omega_n)}{2 \omega_n}. \label{gra}
\end{align} 
\item Update the battery storage threshold via
\begin{align*} 
& B^{(n+1)} = B^{(n)} -\varsigma_n  \partial_{B} J^{(k)}_n
\end{align*} 
which gives
\begin{align*} 
& u^{(n+1)} (P,g,H,B) = \left\{ \begin{array}{cl}
E_0 , & \quad \textrm{if} ~ B \leq B^{(n+1)}  \\
E_1 , & \quad  \textrm{otherwise.}
\end{array} \right. 
\end{align*}
\end{itemize}

The above algorithm is a gradient-estimate 
based algorithm (see \cite{spall2005introduction}) for estimating the 
optimal threshold $B^*$ where only measurements of the loss function is available (i.e., no gradient information). We note that (\ref{gra}) evaluates an approximation to the gradient. This algorithm generates a sequence of estimates for the threshold policy $B^*$ which converges to a local minimum with corresponding energy allocation $u^*$. The reader is referred to \cite{spall2005introduction} for associated convergence analysis of this and other related algorithms (see e.g., Theorem 7.1 in \cite{spall2005introduction}). Note that gradient-estimate based algorithms are sensitive to initial conditions and should be evaluated for several distinct initial conditions to find the best local minimum.

\begin{figure}[!t]
\begin{center}
\includegraphics[width=9cm,height=6cm]{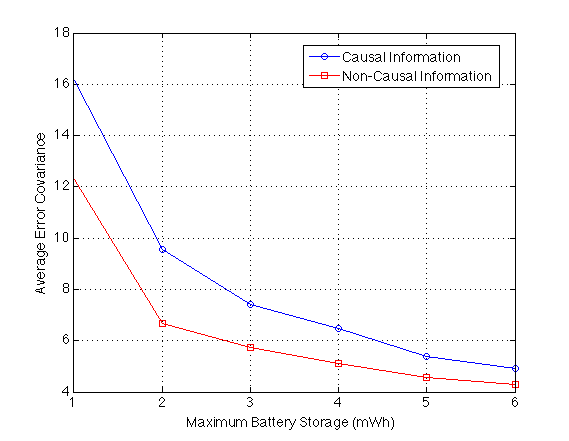}  
\caption{Perfect feedback case: Infinite-time horizon average error covariance versus the maximum battery storage (mWh)} 
\label{fig1:PF-IH-B} 
\end{center}                             
\end{figure}
\begin{figure}[!t]
\begin{center}
\includegraphics[width=9cm,height=6cm]{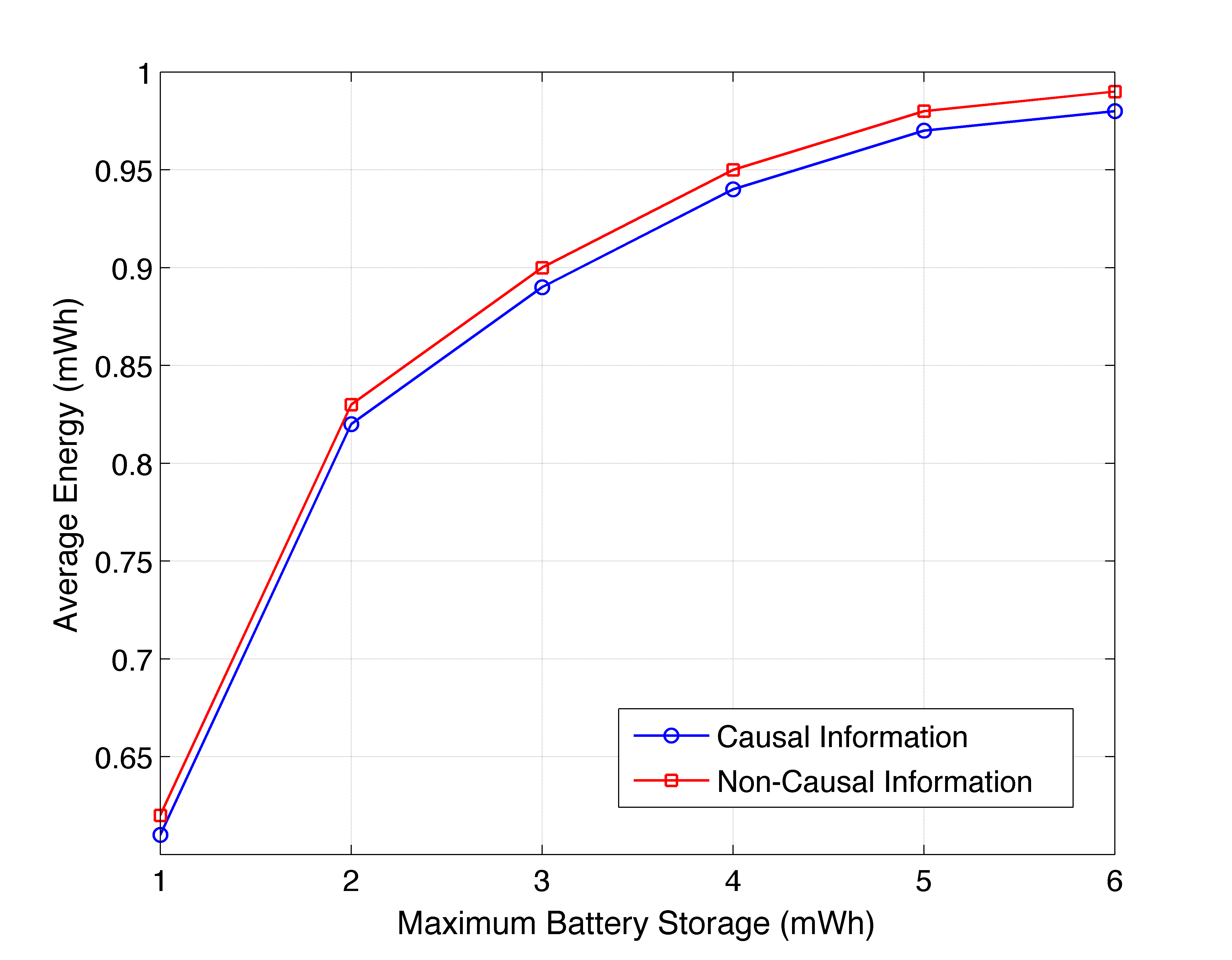}   
\caption{Perfect feedback case: Infinite-time horizon average energy versus the maximum battery storage (mWh)} 
\label{fig2:PF-IH-B-U} 
\end{center}                                
\end{figure}

\section{Numerical Examples} \label{numerical_sec} 

We present here numerical results for a scalar process with following parameters $A=1.2$, $C=1$, $Q=1$, $R=1$ and $P_{x_0}=1$ defined in Section \ref{ProcessDyn_sec}.

In this model we assume that the sensor uses a binary phase shift keying (BPSK) transmission scheme \cite{Proakis} with $b$ bits per packet. Therefore, (\ref{h_defn}) in Section \ref{FCchannel_sec} is of the from
\begin{align*}
& \Prob(\gamma_k = 1|g_k,u_k) = h(g_k u_k) = \Big( \int\limits_{-\infty}^{\sqrt{g_k u_k}} \frac{1}{\sqrt{2 \pi}} e^{-t^2/2} dt \Big)^b
\end{align*}
where we use $b=4$ in the simulations. This model for the packet loss probabilities is studied in \cite{Quevedo_Automatica}.

For simplicity, the fading channel is taken to be Rayleigh \cite{Rappaport} so that $\{g_k\}$ is i.i.d. exponentially distributed with probability density function (p.d.f) of the form $\Prob(g_k) = \frac{1}{\bar{g}} \exp(-g_k/\bar{g})$ with $\bar{g}$ being its mean. We also assume that the harvested energy process $\{H_k\}$ is i.i.d. and exponentially distributed with p.d.f. $\Prob(H_k) = \frac{1}{\bar{H}} \exp(-H_k/\bar{H})$ with $\bar{H}$ being its mean. 

For the following simulation results we use 50 discretization points for each of the quantities of Bellman equations.

\begin{figure}[!t]
\begin{center}
\includegraphics[width=9cm,height=6cm]{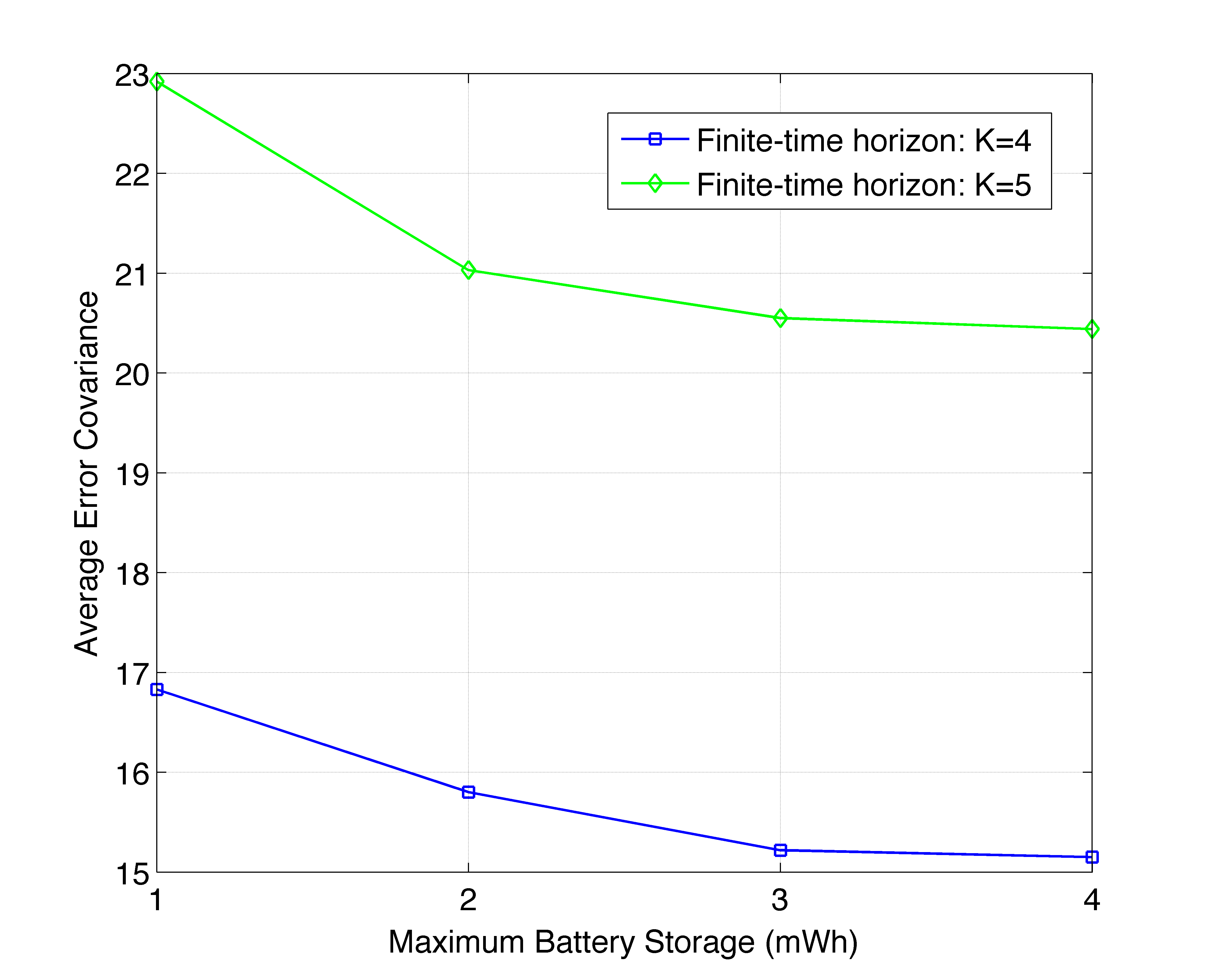}  
\caption{Perfect feedback case: The sum of finite-time horizon expected error covariance versus the maximum battery storage (mWh)} \label{fig8:FTH} 
\end{center}                                
\end{figure}
\begin{figure}[!t]
\begin{center}
\includegraphics[width=9cm,height=6cm]{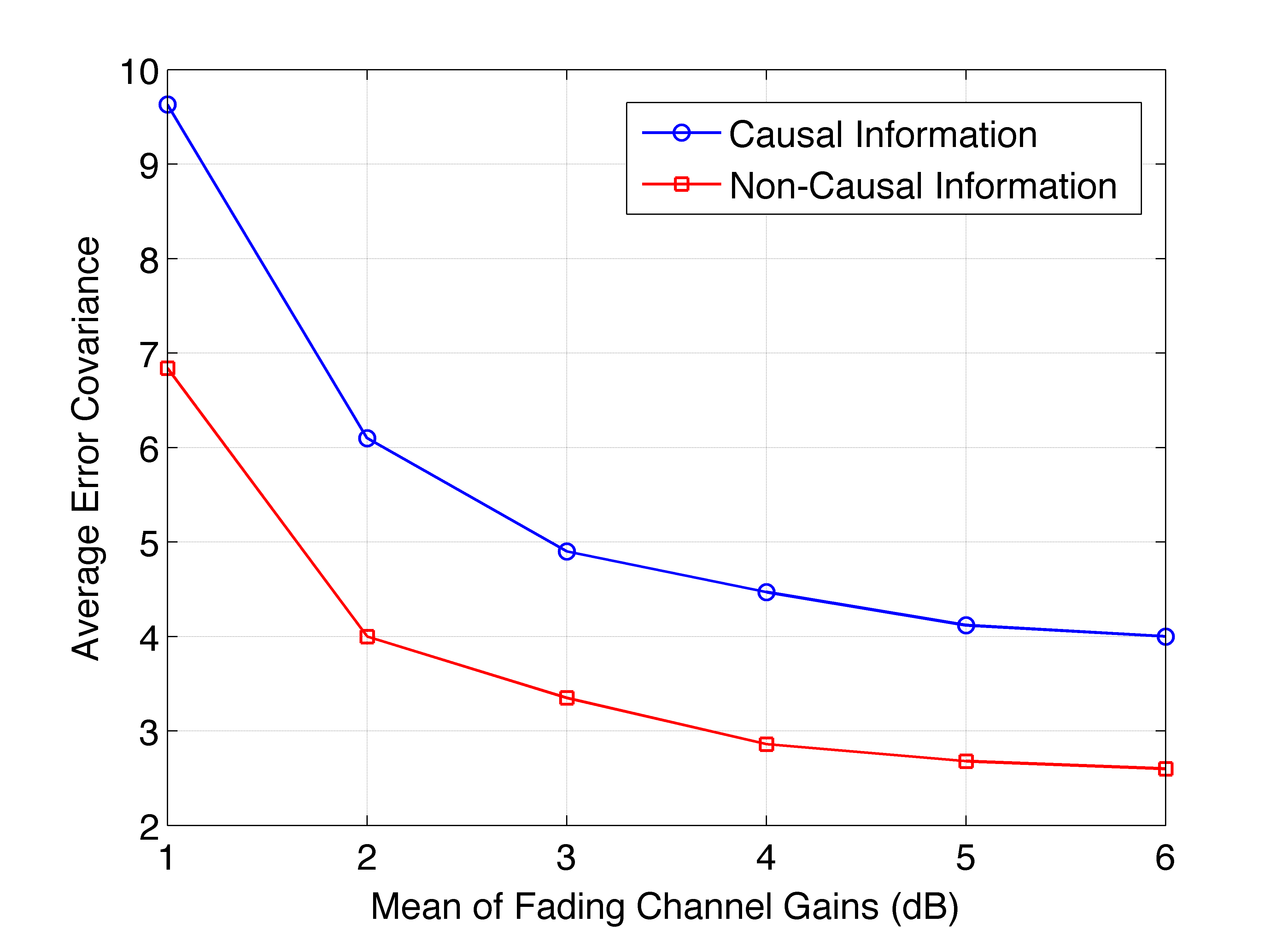}
\caption{Perfect feedback case: Infinite-time horizon average error covariance versus the mean of fading channel gains}
\label{fig3:PF-IH-SNR}
\end{center}
\end{figure}

\subsection{Simulation Results for the Perfect Feedback Communication Channel Case}
 
We first consider the case of perfect packet receipt acknowledgments by setting $\eta$ and $\epsilon$ in Section \ref{ternary_subsec} to zero. We first fix the mean of the fading channel gains to $\bar g = 1$ decibel (dB) and the mean of harvested energy to $\bar H= 1$ milliwatt hour (mWh). Then, we plot the average expected error covariance versus the maximum battery storage energy for the infinite-time horizon formulation (\ref{MDP_constrained:IH}) in Fig. \ref{fig1:PF-IH-B} where both cases of causal and non-causal fading channel gains and energy harvesting information are shown. We note that the performance gets better as the maximum battery storage energy increases in both case. Fig. \ref{fig1:PF-IH-B} also shows that in this setting the performance for the non-causal information case is generally better than the performance of system with causal information.

In Fig. \ref{fig2:PF-IH-B-U} the corresponding average transmission energy versus the maximum battery storage is shown for the infinite-time horizon formulation (\ref{MDP_constrained:IH}) of both cases of causal and non-causal fading channel gains and energy harvesting information. The reader is also referred to Fig. 2 in \cite{Alex_CDC12} that shows the average transmission power versus the expected error covariance trade-off in the case of an average transmission power constraint instead of the energy harvesting constraint considered here. 

For the finite-time horizon case (\ref{MDP_constrained}) when $K=4$ and $5$ the sum of expected error covariance versus the maximum battery storage energy is shown in Fig. \ref{fig1:PF-IH-B} in the case of causal information. As expected the sum of expected error covariances increases when $K$ increases. Similar to graphs of Fig. \ref{fig1:PF-IH-B} the performance gets better as the maximum battery storage energy increases.

We now fix the mean of harvested energy to $\bar H= 1$ (mWh) and the maximum battery storage energy to 2 (mWh). For the infinite-time horizon formulation (\ref{MDP_constrained:IH}) the average expected error covariance versus the mean of the fading channel gains is plotted in Fig \ref{fig3:PF-IH-SNR} for both cases of causal and non-causal information. As shown in Fig \ref{fig3:PF-IH-SNR} the performance gets better as the mean of the fading channel gains increases in both cases.

\begin{figure}[!t]
\centering 
\includegraphics[width=9cm,height=6cm]{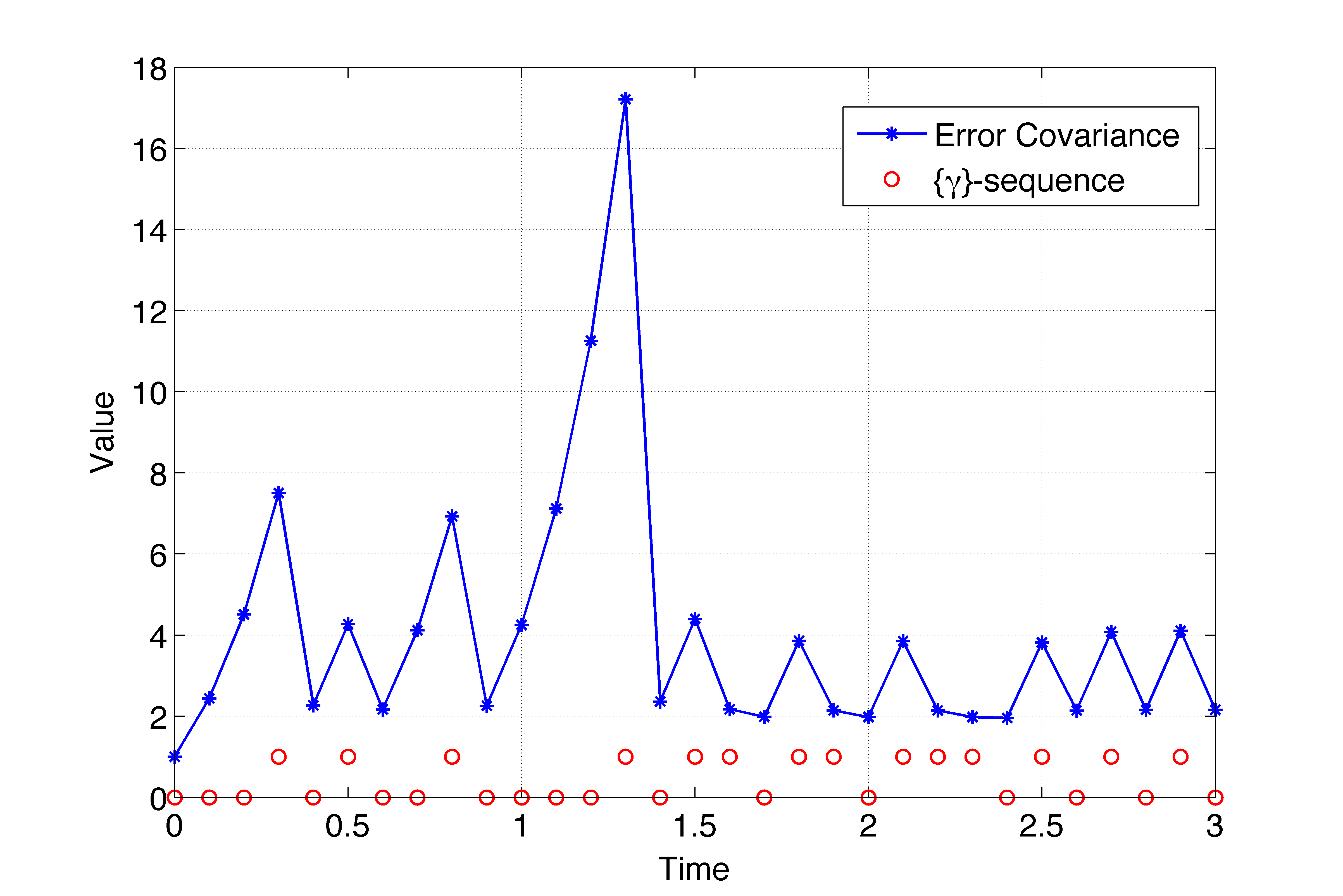} \\
\includegraphics[width=9cm,height=6cm]{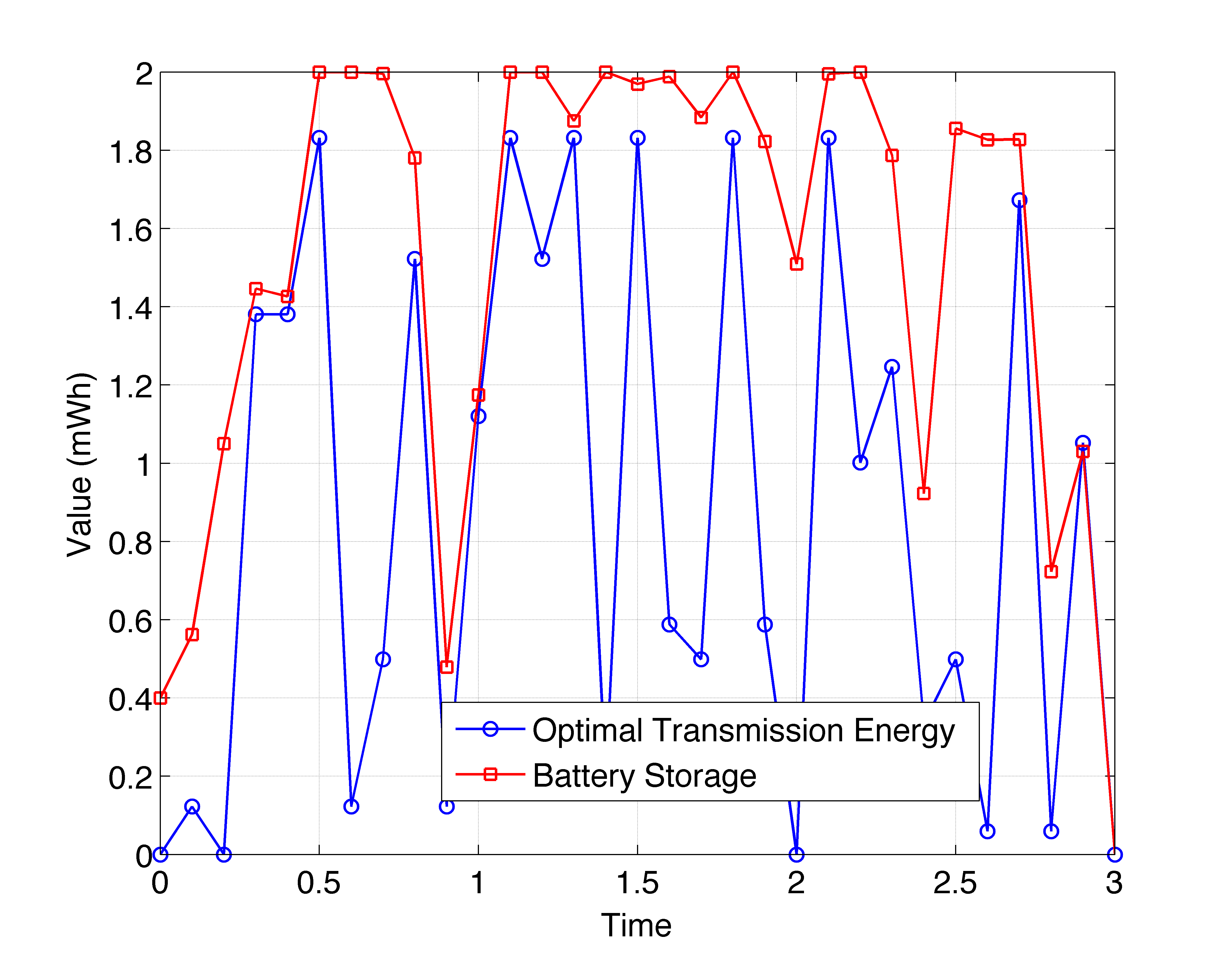}
\caption{Optimal energy allocations}
\label{Fig6-7:power_allocations}
\end{figure}

In Fig. \ref{Fig6-7:power_allocations} we further plot a single simulation run of $\{P_k\}$ with the packet loss process $\{\gamma_k\}$ where $\bar H= 1$ (mWh), $\bar g = 1$ (dB) and $B_{\textrm{max}}=2$ (mWh). The battery storage $\{B_k\}$ and corresponding optimal energy allocations $\{u_k\}$ are also shown.

We can see that in the optimal energy allocation scheme, the allocated energy values will depend not only on the current channel gain $g$ and error covariance $P$ but also on the energy harvesting $H$ through the battery storage $B$. The allocated energy $u$ tends to be higher when the error covariance $P$ is larger, provided the corresponding channel gain $h$ and the battery storage $B$ are not too small.

\subsubsection*{Threshold Policy} 
\begin{figure}[!t]
\begin{center}
\includegraphics[width=9cm,height=6cm]{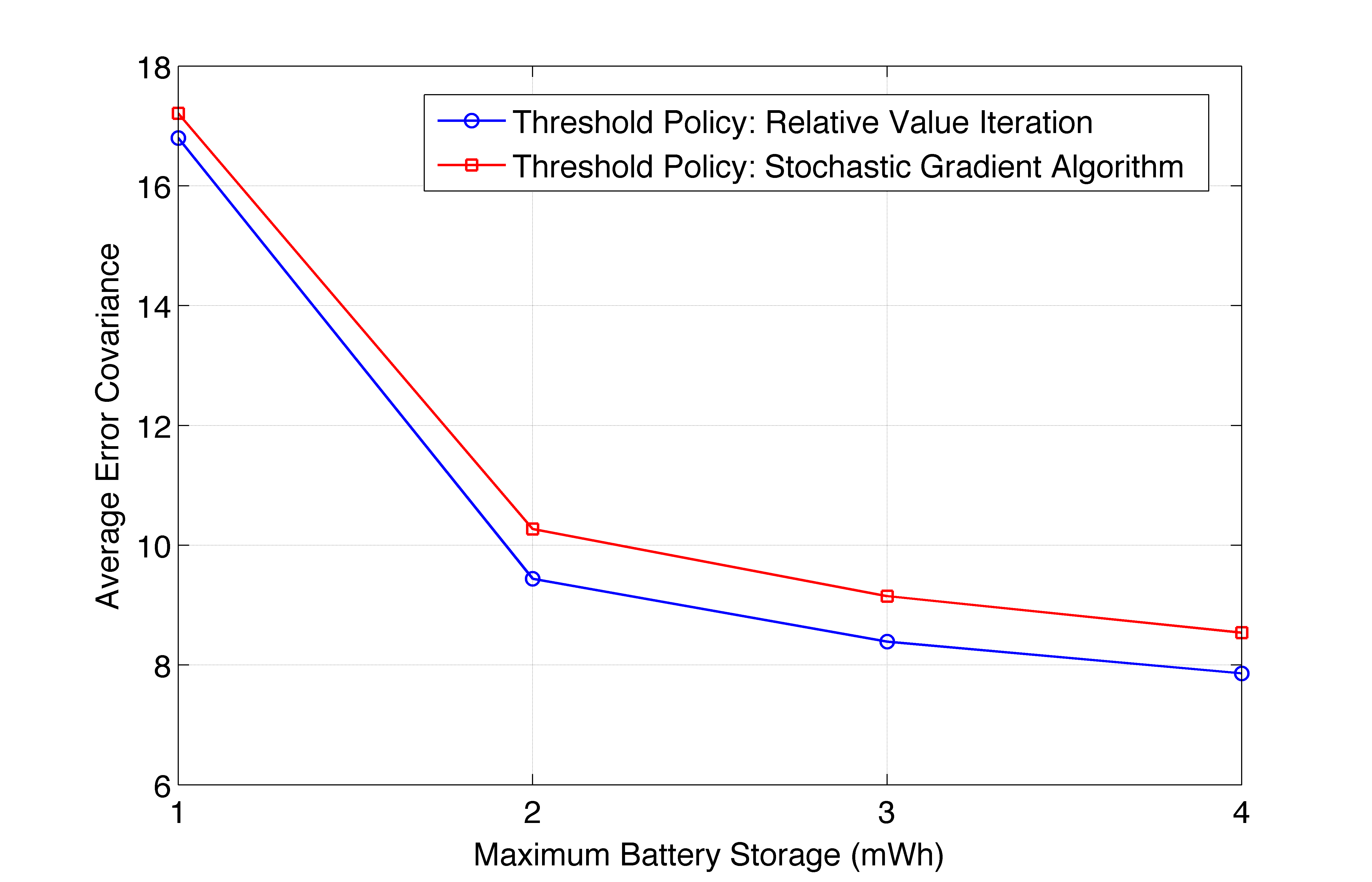}
\caption{Performance of threshold policy using: (i) relative value iteration algorithm, and (ii) stochastic gradient algorithm.}
\label{fig10:Threshold}
\end{center}
\end{figure}

We now consider the case that the transmission energy allocation control belongs to a two element set $\{0,1\}$ in the infinite-time horizon formulation (\ref{MDP_constrained:IH}) with parameters $\bar H= 1$ (mWh), $\bar g = 1$ (dB) and $B_{\textrm{max}}=2$ (mWh). As explained at the end of Section \ref{Str_sec} the optimal transmission energy allocation policy is threshold of the form 
\begin{align*} 
& u^* (P,g,H,B) = \left\{ \begin{array}{cl}
0 , & \quad \textrm{if} ~ B \leq B^*(P,g,H,B)   \\
1 , & \quad  \textrm{otherwise}
\end{array} \right. 
\end{align*}
where $B^*(\cdot,\cdot,\cdot,\cdot)$ is the corresponding battery storage threshold. Applying the stochastic gradient algorithm of Section \ref{Str_sec} with parameters $\omega=0.1, \varsigma =0.5$ and $\kappa=1$ to our model yields a set of threshold policies which gives $u^*$. 
Fig. \ref{fig10:Threshold} shows the simulation results where the relative value iteration algorithm and the threshold policy based algorithm are used. It can be seen that there is a small gap between the 
simulation results obtained via the two methods. This can be attributed to the fact that the optimal threshold is not exactly calculated by the stochastic gradient algorithm which only converges to a local minimum. 

\begin{figure}[!t]
\begin{center}
\includegraphics[width=9cm,height=6cm]{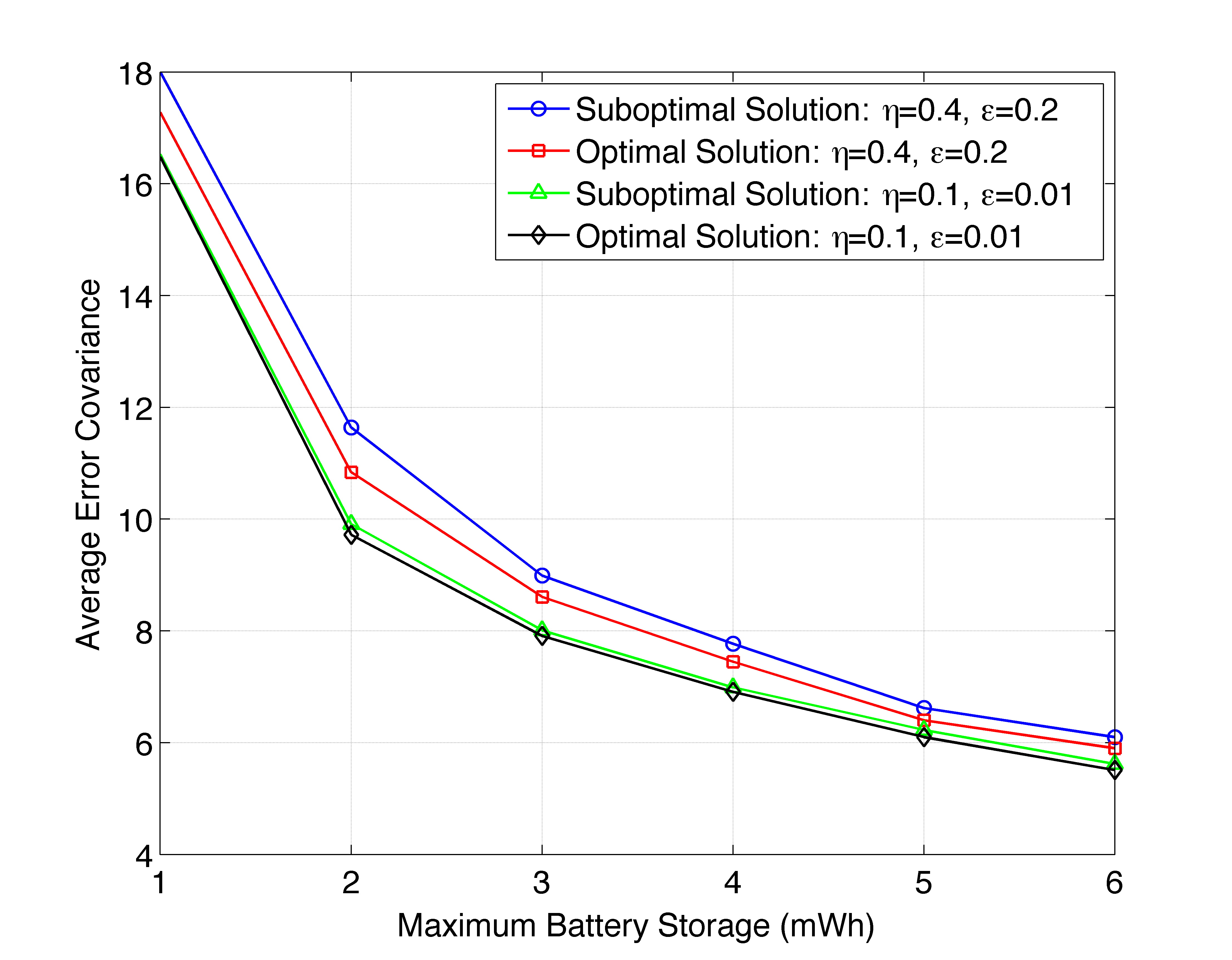}  
\caption{Imperfect feedback case (optimal and suboptimal solutions for parameters: (i) $\eta=0.4$ and $\epsilon=0.2$, and (ii) $\eta=0.1$ and $\epsilon=0.01$): Infinite-time horizon average error covariance versus the maximum battery storage (mWh)} 
\label{Fig4:IF-IH-B} 
\end{center}                             
\end{figure}

\subsection{Simulation Results for the Imperfect Feedback Communication Channel Case}

We now consider the case of imperfect packet receipt acknowledgments as given in Section \ref{ternary_subsec} with two sets of parameters: (i) $\eta=0.4$ and $\epsilon=0.2$, and (ii) $\eta=0.1$ and $\epsilon=0.01$. In this simulation we fix the mean of the fading channel gains to $\bar g = 1$ (dB) and the mean of harvested energy to $\bar H= 1$ (mWh). Then, we plot in Fig. \ref{Fig4:IF-IH-B} both optimal and suboptimal results of the average expected error covariance versus the maximum battery storage energy for the infinite-time horizon formulation (\ref{MDP_constrained:IH}). Similar to Fig. \ref{fig1:PF-IH-B} the performance gets better as the maximum battery storage energy increases. Fig. \ref{Fig4:IF-IH-B} also shows that, as expected, the performance for the optimal solution (see Section \ref{optimSol_sec}) is generally better than the performance of the suboptimal solution (see Section \ref{Sub_sec}). But, for small packet loss and error probabilities of the feedback communication channel, $\eta$ and $\epsilon$, the performance for the suboptimal solution, which is easier to implement, is close to the performance of optimal solution. Note that optimal solutions are computationally demanding since it is required to solve Bellman equations on a discretized subset of the space of probability densities.

\section{Conclusions} \label{conc_sec} We have studied the problem of optimal transmission energy allocation for estimation error covariance minimization in Kalman filtering with random packet losses over a fading channel when the sensor is equipped with energy harvesting technology. The feedback channel from receiver to sensor is an erroneous packet dropping link which models an imperfect receipt acknowledgments channel. In this problem formulation either a sum over a finite-time horizon or a long term average (infinite-time horizon) of the trace of the expected estimation error covariance of the Kalman filter is minimized, subject to energy harvesting constraints. The resulting Markov decision process problems with imperfect state information are solved by the use of the dynamic programming principle. Using the concept of submodularity, some structural results on the optimal transmission energy allocation policy are obtained. Suboptimal solutions are also discussed which are far less computationally intensive than optimal solutions.

\begin{appendix}
\renewcommand{\theequation}{A.\arabic{equation}}
\setcounter{equation}{0}
{\textit{Proof of Lemma \ref{lem:information-state}:}} The total probability
formula\footnote{$ \Prob(A\,|\,B)=\sum_i  \Prob(A,C_i\,|\,B)$} and the chain rule give  
\begin{align}
& \Prob(P_{k+1},z^k,g_k,u_k) \! = \! \sum_{\gamma_k} \! \int_{P_k} \!\! \Prob (P_{k+1},P_k,\gamma_k,z^k,g_k,u_k) dP_k \notag \\
& = \!  \sum_{\gamma_k} \! \int_{P_k} \!\!\! \Prob (P_{k+1}|P_k,\gamma_k,z^k,g_k,u_k) \Prob (P_k,\gamma_k,z^k,g_k,u_k)  dP_k \notag \\
& = \! \sum_{\gamma_k} \! \int_{P_k} \!\!\! \Prob (P_{k+1}|P_k,\gamma_k) \Prob (P_k,\gamma_k,z^k,g_k,u_k)  dP_k \label{first_equalit}
\end{align}
where the last equality is because $P_{k+1}$ is a function of $P_k$ and $\gamma_k$ by (\ref{EH:RiccatiEq}). But, the chain rule implies that
\begin{align}
& \Prob (P_k,\gamma_k,z^k,g_k,u_k) \! = \! \Prob (P_k,\gamma_k,z^{k-1},\hat \gamma_k, g_{k-1},u_{k-1},g_k,u_k) \notag \\
& \hspace{0.1cm} = \Prob (\hat \gamma_k | P_k,\gamma_k,z^{k-1},g_{k-1},u_{k-1},g_k,u_k) \notag \\
& \hspace{0.5cm} \times \Prob (\gamma_k | P_k,z^{k-1},g_{k-1},u_{k-1},g_k,u_k) \notag \\
& \hspace{0.5cm} \times \Prob (g_k | P_k,z^{k-1},g_{k-1},u_{k-1},u_k) \notag \\
&  \hspace{0.5cm} \times  \Prob (P_k|z^{k-1},g_{k-1},u_{k-1},u_k) \Prob (z^{k-1},g_{k-1},u_{k-1},u_k) \notag \\
& \hspace{0.1cm} = \Prob (\hat \gamma_k | \gamma_k) \Prob (\gamma_k |g_k,u_k) \Prob (g_k |g_{k-1}) \Prob (P_k|z^{k-1},g_{k-1},u_{k-1}) \notag \\
&  \hspace{0.5cm} \times \Prob (z^{k-1},g_{k-1},u_{k-1},u_k).  \label{second_equalit}
\end{align}

Substituting (\ref{second_equalit}) in (\ref{first_equalit}) yields
\begin{align}
& \Prob(P_{k+1},z^k,g_k,u_k) = \sum_{\gamma_k} \int_{P_k}  \Big(\Prob (P_{k+1}|P_k,\gamma_k) \Prob (\hat \gamma_k | \gamma_k) \notag \\
& \hspace{0.2cm} \times \Prob (\gamma_k |g_k,u_k) \Prob (g_k |g_{k-1}) \Prob (P_k|z^{k-1},g_{k-1},u_{k-1}) \notag \\
& \hspace{0.2cm} \times  \Prob (z^{k-1},g_{k-1},u_{k-1},u_k) \Big) dP_k. \label{simcommon}
\end{align}

On the other hand, 
\begin{align}
& \Prob(P_{k+1}|z^k,g_k,u_k) = \alpha \times \Prob(P_{k+1},z^k,g_k,u_k)  \label{alphaval}
\end{align}
where $\alpha$ is a normalizing constant. Integrating (\ref{alphaval}) with respect to $P_{k+1}$ 
\linebreak
gives 
$\alpha = \big(\int_{P_{k+1}} \Prob(P_{k+1},z^k,g_k,u_k) d P_{k+1}\big)^{-1}.$ But,
\begin{align}
& \int_{P_{k+1}} \Prob(P_{k+1},z^k,g_k,u_k) d P_{k+1} \notag \\
&  =   \int_{P_{k+1}} \Big[\sum_{\gamma_k} \int_{P_k}  \Big(\Prob (P_{k+1}|P_k,\gamma_k) \Prob (\hat \gamma_k | \gamma_k)\Prob (\gamma_k |g_k,u_k) \notag \\
& \hspace{1cm} \times \Prob (g_k |g_{k-1}) \Prob (P_k|z^{k-1},g_{k-1},u_{k-1}) \notag \\
& \hspace{1cm}  \times \Prob (z^{k-1},g_{k-1},u_{k-1},u_k) \Big) dP_k \Big] dP_{k+1}. \label{normval}
\end{align}
By changing the order of integration, we may simplify (\ref{normval}) as
\begin{align}
& \int_{P_{k+1}} \Prob(P_{k+1},z^k,g_k,u_k) d P_{k+1} \notag \\
& \hspace{0.2cm} =   \Prob (g_k |g_{k-1}) \Prob (z^{k-1},g_{k-1},u_{k-1},u_k) \notag \\
& \hspace{0.5cm}   \times \sum_{\gamma_k} \int_{P_k}  \Big( \big( \int_{P_{k+1}} \Prob (P_{k+1}|P_k,\gamma_k) dP_{k+1} \big) \Prob (\hat \gamma_k | \gamma_k) \notag \\
&   \hspace{2cm} \times \Prob (\gamma_k |g_k,u_k)\Prob (P_k|z^{k-1},g_{k-1},u_{k-1})  \Big) dP_k \notag \\
& \hspace{0.2cm} =   \Prob (g_k |g_{k-1}) \Prob (z^{k-1},g_{k-1},u_{k-1},u_k) \sum_{\gamma_k} \Big( \Prob (\hat \gamma_k | \gamma_k)\notag \\
& \hspace{0.5cm}   \times \Prob (\gamma_k |g_k,u_k) \big( \int_{P_k}\Prob (P_k|z^{k-1},g_{k-1},u_{k-1}) dP_k \big) \Big) \notag \\
& \hspace{0.2cm} =   \Prob (g_k |g_{k-1}) \Prob (z^{k-1},g_{k-1},u_{k-1},u_k) \notag \\
& \hspace{1cm} \times \sum_{\gamma_k} \Prob (\hat \gamma_k | \gamma_k)\Prob (\gamma_k |g_k,u_k)
\end{align}
where we used the fact that $\int_{P_{k+1}} \Prob (P_{k+1}|P_k,\gamma_k) dP_{k+1} =1$ and $\int_{P_k}\Prob (P_k|z^{k-1},g_{k-1},u_{k-1}) dP_k=1$. Hence, we have
\begin{align}
&\alpha = \Big(\Prob (g_k |g_{k-1}) \Prob (z^{k-1},g_{k-1},u_{k-1},u_k) \notag \\
& \hspace{2cm} \times \sum_{\gamma_k} \Prob (\hat \gamma_k | \gamma_k)\Prob (\gamma_k |g_k,u_k)\Big)^{-1}.  \label{alphasim}
\end{align}

Finally, substituting (\ref{simcommon}) and (\ref{alphasim}) in (\ref{alphaval}) gives
\begin{align*}
& \Prob(P_{k+1}|z^k,g_k,u_k) \\
& = \sum_{\gamma_k} \Big[ \int_{P_k}  \!\!\Big(\Prob (P_{k+1}|P_k,\gamma_k) \times \Prob (P_k|z^{k-1},g_{k-1},u_{k-1}) \Big) dP_k \\
& \hspace{1cm} \times \frac{\Prob (\hat \gamma_k | \gamma_k)\Prob (\gamma_k |g_k,u_k)}{\sum_{\gamma_k} \Prob (\hat \gamma_k | \gamma_k)\Prob (\gamma_k |g_k,u_k)} \Big]
\end{align*}
as the information-state recursion given in (\ref{InfState_recur_tot}). \qed

{\it Proof of Theorem \ref{Causal:IH:Bellman}}: 
We first show the inequality
\begin{align}
& \rho + V(\pi,g,H,B) \geq \min_{0 \leq u \leq B} \Big\{\mathbb{E} \big[\tr\big(\mathcal L(P,\gamma)\big)\big|\pi,g,u\big] \notag \\
& \hspace{0.5cm} + \mathbb{E} \Big[V\big(\Phi \big[\pi,\hat \gamma,g,u\big],\tilde g,\tilde H,\min\{B-u+\tilde H, B_{\textrm{max}}\}\big) \notag \\
& \hspace{3cm} \big| \pi,g,H,u\Big] \Big\} \label{ACOI}
\end{align}
by verifying conditions (W) and (B) of  \cite{Schal} that guarantee the existence of solutions to (\ref{ACOI}) for MDPs with general state space. Denote the state space $\mathcal{S}$ and action space $\mathcal{A}$, i.e. $(\pi_k,g_k,H_k,B_k) \in \mathcal{S}$ and  $u_k \in \mathcal{A}$. 
Condition (W) of \cite{Schal} in our notation says that:\\
0) The state space $\mathcal{S}$ is locally compact.\\
1) Let $U(\cdot)$ be the mapping that assigns to each $(\pi_k,g_k,H_k,B_k)$ the nonempty set of available actions. Then $U(\pi_k,g_k,H_k,B_k)$ lies in a compact subset of $\mathcal{A}$ and $U(\cdot)$ is upper semicontinuous.\\
2) The transition probabilities are weakly continuous.\\\
3) $\mathbb{E} \big[\tr\big(\mathcal L(P,\gamma)\big)\big|\pi,g,u\big]$ is lower semicontinuous. \\
By our assumption that $u_k \leq B_k \leq B_{max}$, 0) and 1) of (W) can be easily verified. The condition 2) follows from (\ref{FunPhi}), while condition 3) follows from the definition (\ref{EH:RiccatiEqOperator}). 

We define $w_{\delta}(\pi_0,g_0,H_0,B_0) = v_{\delta} (\pi_0,g_0,H_0,B_0) - m_{\delta}$ where
\begin{align*}
& v_{\delta} (\pi_0,g_0,H_0,B_0) = \inf_{\{u_k: k \geq 0\}} \mathbb{E}[\sum_{k=0}^\infty \delta^k \mathbb{E} \big[\tr\big(\mathcal L(P_k,\gamma_k)\big) \\
& \hspace{4.5cm} \big|\pi_k,g_k,u_k\big]|\pi_0,g_0,H_0,B_0]
\end{align*}
and $m_{\delta} = \inf_{(\pi_0,g_0,H_0,B_0)} v_\delta (\pi_0,g_0,H_0,B_0)$, then Condition (B) of \cite{Schal} in our notation implies that 
\begin{align*}
& \sup_{\delta<1} w_{\delta} (\pi_0,g_0,H_0,B_0) < \infty, \qquad \forall ~ (\pi_0,g_0,H_0,B_0).
\end{align*}

Following Section 4 of \cite{Schal}, we define the stopping time 
$\tau = \inf \{k \geq 0: v_{\delta} (\pi_k,g_k,H_k,B_k) \leq m_{\delta} + \varsigma\}$ for some $\varsigma \geq 0$. 
Given $\varsigma > 0$ and an arbitrary $(\pi_0,g_0,H_0,B_0)$, consider a suboptimal power allocation policy where the sensor transmits based on the same policy as the one that achieves $m_\delta$ (with a different initial condition) until $v_{\delta} (\pi_N,g_N,H_N,B_N) \leq m_{\delta} + \varsigma$ is satisfied at some time $N$. By the exponential forgetting property of initial conditions for Kalman filtering, we have $N < \infty$ with probability $ 1$ and $\mathbb{E}[N] < \infty$. Since $\tau \leq N$, we have $\mathbb{E}[\tau] < \infty$. Then by Lemma 4.1 of \cite{Schal},
\begin{align}
& w_\delta (\pi_0,g_0,H_0,B_0) \leq \varsigma   \notag \\
& \quad + \inf_{\{\gamma_k\}} \mathbb{E}[\sum_{k=0}^{\tau-1} \mathbb{E} \big[\tr\big(\mathcal L(P_k,\gamma_k)\big) \big|\pi_k,g_k,u_k\big]|\pi_0,g_0,H_0,B_0]  \notag \\ 
& \leq \varsigma + \mathbb{E}[\tau] \times Z  < \infty \label{w_beta_bound}
\end{align}
where the second inequality uses Wald's equation, with $Z$ being an upper to the expected error covariance that exists by Theorem \ref{Stability:thm}. Hence condition (B) of \cite{Schal} is satisfied and a solution to (\ref{ACOI}) exists.  

To show equality in (\ref{ACOI}), we will require a further equicontinuity property to be satisfied. This can be shown by a similar argument as in the proof of Proposition 3.2 of \cite{HuangDey_TSP}. The assumptions in Sections 5.4 and 5.5 of \cite{HernandezLermaLasserre} may then be verified to conclude the existence of a solution to the average cost optimality equation (\ref{perfecfeed_dp:IH}). 

\end{appendix}

\bibliographystyle{IEEEtran}
\bibliography{FP-13-397-Ver2-ref}   

\end{document}